\documentclass[11pt]{article}
\usepackage{amsmath}
\usepackage{color}
\newcommand{\Limsup}{\mathop{{\rm Lim}\,{\rm sup}}}
\def\disp{\displaystyle}

\def\tto{\;{\lower 1pt \hbox{$\rightarrow$}}\kern-10pt
\hbox{\raise 2pt \hbox{$\rightarrow$}}\;}
\def\Hat{\widehat}
\def\hat{\widehat}
\def\Tilde{\widetilde}

\def\Bar{\overline}
\def\ra{\rangle}
\def\la{\langle}
\def\ve{\varepsilon}
\def\B{I\!\!B}
\def\IN{I\!\!N}

\def\h{\hfill\Box}
\def\R{I\!\!R}

\def\K{\mathcal{K}}
\def\ox{\bar{x}}
\def\op{\bar{p}}

\def\oy{\bar{y}}
\def\oz{\bar{z}}
\def\ov{\bar{v}}

\def\ow{\bar{w}}

\def\inte{\mbox{\rm int}\,}

\def\gph{\mbox{\rm gph}\,}

\def\epi{\mbox{\rm epi}\,}
\def\dim{\mbox{\rm dim}\,}
\def\dom{\mbox{\rm dom}\,}

\def\h{\hfill\triangle}
\def\dn{\downarrow}
\def\O{\Omega}
\def\Lm{\Lambda}
\def\ph{\varphi}
\def\emp{\emptyset}
\def\st{\stackrel}
\def\oR{\Bar{\R}}

\def\lm{\lambda}
\def\gg{\gamma}
\def\dd{\delta}
\def\al{\alpha}
\def\kk{\kappa}

\def\Th{\Theta}
\def\th{\theta}
\def\vt{\vartheta}
\newcounter{lk}
\def\Limsup{\mathop{{\rm Lim}\,{\rm sup}}}

\begin{document}
\newtheorem{Theorem}{Theorem}[section]
\newtheorem{Conjecture}[Theorem]{Conjecture}
\newtheorem{Proposition}[Theorem]{Proposition}
\newtheorem{Remark}[Theorem]{Remark}
\newtheorem{Lemma}[Theorem]{Lemma}
\newtheorem{Corollary}[Theorem]{Corollary}
\newtheorem{Definition}[Theorem]{Definition}
\newtheorem{Example}[Theorem]{Example}
\renewcommand{\theequation}{\thesection.\arabic{equation}}
\def\proof{
\normalfont
\medskip
{\noindent\itshape Proof.\hspace*{6pt}\ignorespaces}}
\def\endproof{$\h$ \vspace*{0.1in}}
\begin{center}
\vspace*{0.3in} {\bf LOCAL STRONG MAXIMAL MONOTONICITY\\ AND FULL STABILITY FOR PARAMETRIC VARIATIONAL SYSTEMS} \\[2ex]
B. S. MORDUKHOVICH\footnote{Department of Mathematics, Wayne State University, Detroit, MI 48202 (boris@math.wayne.edu). Research of this author was partially supported by the National Science Foundation under grant DMS-1007132,}
and T. T. A. NGHIA\footnote{Department of Mathematics and Statistics, Oakland University, Rochester, MI 48309 (nttran@oakland.edu). Research of this author was partially supported by postdoctoral fellowships from the Pacific Institute for the Mathematical Sciences and the University of British Columbia.}
\end{center}
\small{\sc Abstract.} The paper introduces and characterizes new notions of Lipschitzian and H\"olderian full stability of solutions to general parametric variational systems described via partial subdifferential and normal cone mappings acting in Hilbert spaces. These notions, postulated certain quantitative properties of single-valued localizations of solution maps, are closely related to local strong maximal monotonicity of associated set-valued mappings. Based on advanced tools of variational analysis and generalized differentiation, we derive verifiable characterizations of the local strong maximal monotonicity and full stability notions under consideration via some positive-definiteness conditions involving second-order constructions of variational analysis. The general results obtained are specified for important classes of variational inequalities and variational conditions in both finite and infinite dimensions.\\[1ex]
{\em 2010 Mathematics Subject Classification}. Primary 49J53; Secondary 49J52, 90C31.\\[1ex]
{\em Key words and phrases}. Variational analysis, parametric variational systems, variational inequalities and variational conditions, local monotonicity and strong maximal monotonicity, full stability, Legendre forms, polyhedricity, generalized differentiation, subdifferentials, coderivatives.

\normalsize
\section{Introduction}
\setcounter{equation}{0}

The paper belongs to the area of modern variational analysis, which has been well recognized as a fruitful field of mathematics with numerous applications; see, e.g., the books \cite{DR,M1,rw} and the references therein. We pursue here a twofold goal: to study {\em local strong maximal monotonicity} of set-valued operators in Hilbert spaces and the usage of ideas and results on strong  maximal monotonicity to introduce and characterize new notions of {\em quantitative stability}, in both Lipschitzian and H\"olderian settings, for {\em parametric variational systems} (PVS) given in the general form
\begin{equation}\label{1.1}
v\in f(x,p)+\partial_x g(x,p)
\end{equation}
as well as its important specifications. In \eqref{1.1} we have: $x\in X$ is the {\em decision} variable from a Hilbert space $X$; $(v,p)\in X\times P$ is a pair of perturbation {\em parameters}, where $v\in X$ signifies {\em canonical} perturbations while $p\in P$ stands for {\em basic} perturbations taking values in a metric space $P$; the single-valued {\em base} mapping $f\colon X\times P\to X$ is smooth around the reference point $(\ox,\op)$; and the {\em potential} $g\colon X\times P\to\oR:=(-\infty,\infty]$ is an extended-real-valued and lower semicontinuous (l.s.c.) function with the symbol $\partial_x$ indicating the set of its {\em partial limiting subgradients} with respect to the decision variable; see Section~2 for more details.

It has been realized over the years that model \eqref{1.1} and its various specifications (known, in particular, as variational and quasi-variational inequalities, generalized equations, and variational conditions) provide convenient frameworks for the study and applications of many important issues of nonlinear analysis, partial differential equations, optimization, equilibria, control theory, numerical algorithms, etc.; see, e.g., \cite{DR,FP,hms,KS,M1,R0,R03,rw} and the bibliographies therein as well as further references presented below in this paper.

The vast majority of research on parametric variational systems revolves around establishing certain stability properties of the {\em solution map}
\begin{equation}\label{1.2}
S(v,p):=\big\{x\in X\big|\;v\in f(x,p)+\partial_x g(x,p)\big\},\quad(v,p)\in X\times P,
\end{equation}
to PVS \eqref{1.1} with respect to perturbations of the reference parameter pair $(\ov,\op)$. Starting with the pioneering papers by Stampacchia \cite{S} for variational inequalities (motivated by applications to partial differential equations) and by Robinson \cite{R0} for generalized equations (motivated by applications to optimization), many publications in this direction have been devoted to deriving efficient conditions ensuring the {\em single-valuedness} and {\em continuity} or {\em Lipschitz continuity} of solution maps \eqref{1.2} to important specifications of PVS with their further applications to various fields of mathematics including those mentioned above. It is worth mentioning that the main tools of analysis in the aforementioned developments were related to the usage of fundamental results from implicit function and topological degree theories; see, e.g., the books \cite{DR,FP} and their references.

Another approach to Lipschitzian stability of solution maps to PVS \eqref{1.1} and more general types of parameter-dependent generalized equations was initiated by the first author \cite{M0} who employed the machinery of nonsmooth variational analysis based on his {\em coderivative characterization} \cite{m93} of the {\em Lipschitz-like} (or Aubin's, pseudo-Lipschitz) property of multifunctions together with well-developed coderivative calculus. However, the main results of \cite{M0} were concerned Lipschitzian stability of {\em set-valued} solution maps while their {\em single-valuedness} was established therein only by imposing rather restrictive {\em monotonicity assumptions} on the initial data of the generalized equations under consideration. On the other hand, it has been shown by Dontchev and Rockafellar \cite{DR1} that, under some smoothness requirements on $f$, the Lipschitz-like property is {\em equivalent} to the simultaneous validity of the local single-valuedness and the classical Lipschitz continuity of solution maps to finite-dimensional variational inequalities over {\em polyhedral} convex sets, which correspond to \eqref{1.2} in the case when $g(x)=\dd_C(x)$ is the indicator function of a convex polyhedron $C\subset\R^n$. Furthermore, a certain ``critical face" characterization of these equivalent properties of solution maps to such variational inequalities was established in \cite{DR1} by using the aforementioned coderivative criterion \cite{m93} via a suitable linearization procedure.\vspace*{0.05in}

The major goal of this paper is to employ advanced tools of first-order and second-order variational analysis and generalized differentiation to deriving verifiable {\em characterizations} of new notions of quantitative stability for PVS \eqref{1.1} and their remarkable specifications. These stability notions {\em imply} (being properly stronger than) the local {\em single-valuedness} and {\em Lipschitz} or {\em H\"older continuity} of the solution maps \eqref{1.2} without a priori monotonicity assumptions imposed on the initial data of PVS. When $f(x,p)=0$ in \eqref{1.1}, the obtained characterizations allow us to conclude that the new stability notions for PVS are equivalent to Lipschitzian (resp.\ H\"olderian) full stability for local minimizers of $g$ introduced by Levy, Poliquin and Rockafellar \cite{LPR} (resp.\ by Mordukhovich and Nghia \cite{MN}). Based on this, we label the new stability notions for PVS as ``full Lipschitzian and H\"olderian stability" of the corresponding solution maps, observing then that these concepts are different from the standard local single-valuedness and Lipschitz or H\"older continuity of \eqref{1.2}; see Section~4 for details. Furthermore, it occurs that in the absence of the parameter $p$ in \eqref{1.1} both notions of Lipschitzian and H\"olderian full stability for \eqref{1.1} reduces to the local strong maximal monotonicity of the inverse solution mapping $S^{-1}$ characterized in this paper. It indicates that this kind of {\em local monotonicity} is behind the full quantitative stability notions under consideration.\vspace*{0.02in}

The rest of our paper is organized as follows. In Section~2 we recall some basic notions of variational analysis and generalized differentiation used in the paper. Section~3 is devoted to a systematic study of {\em local strong maximal monotonicity} of set-valued mappings in Hilbert spaces. We establish there several neighborhood and pointwise {\em coderivative characterizations} of this property and discuss some related results. The section is self-contained, and the results obtained therein are of their own interest, while they are very instrumental to proceed further with the study of quantitative full stability of general parametric variational systems and their subsequent specifications.

Section~4 is central in the paper. We introduce and discuss there the notions of {\em H\"olderian} and {\em Lipschitzian full stability} for general PVS \eqref{1.1} and derive complete {\em second-order characterizations} of these stability notions in both the {\em neighborhood form} (using generalized differential constructions in a neighborhood of the reference point) in the case of infinite-dimensional spaces and the {\em  pointwise form} (using only the point in question) when the decision and parameter spaces are finite-dimensional. The characterizations obtained are expressed in terms of certain {\em positive-definiteness} conditions involving appropriate second-order subdifferential constructions.

Section~5 concerns parametric {\em variational inequalities} written in the generalized equation form
\begin{eqnarray}\label{vi}
v\in f(x,p)+N_C(x)\;\mbox{ with }\;x\in C\subset X,\;p\in P,
\end{eqnarray}
where $N_C(x)$ is the normal cone at $x$ to the closed and {\em convex} subset $C$ of the Hilbert space $X$, and where $P$ is a metric space. It is clear that \eqref{vi} is a particular case of PVS \eqref{1.1} with $g(x)=\delta_C(x)$ being the indicator function of the parameter-independent convex set $C$. By definition of the normal cone in convex analysis, \eqref{vi} can be rewritten in the standard form of parameterized variational inequalities: given $p\in P$, find $x\in C$ such that
\begin{eqnarray}\label{vi1}
\la v-f(x,p),u-x\ra\le 0\;\mbox{ for all }\;u\in C.
\end{eqnarray}
Infinite-dimensional variational inequalities in form \eqref{vi1} often appear in optimization-related (in particular, optimal control) problems governed by elliptic partial differential equations, which are usually modeled via the so-called {\em Legendre form} under the {\em polyhedricity} assumption on $C$; see, e.g., \cite{BBS,B,BS,hms,HS,IT,KS} and the precise definitions in Section~5 for more details. Imposing these natural requirements and employing the results of Section~4, we derive in Section~5 {\em pointwise} necessary and sufficient conditions for Lipschitzian full stability of solutions to the perturbed variational inequalities \eqref{vi} held in the infinite-dimensional framework of Hilbert decision spaces.

In Section~6 we study the parametric {\em variational conditions} given by
\begin{eqnarray}\label{qvi}
v\in f(x,p)+N_{C(p)}(x)\;\mbox{ with }\;x\in C(p)\subset X,\;p\in P
\end{eqnarray}
via the limiting normal cone $N_{C(p)}$ to $C(p)$ at $x$ (see Section~2), where both the decision space $X$ and the parameter space $P$ are finite-dimensional, and where the parameter-dependent set $C(p)$ is described by the finitely many inequality constraints
\begin{eqnarray}\label{ine}
C(p):=\big\{x\in X\big|\;\ph_i(x,p)\le 0\;\mbox{ for }\;i=1,\ldots,m\big\}
\end{eqnarray}
defined by ${\cal C}^2$-smooth functions. The parametric variational conditions \eqref{qvi}, known also as generalized equations with parameter-dependent multivalued terms/fields, etc., are imbedded into the PVS framework \eqref{1.1} with $g(x,p)=\dd_{C(p)}(x)$ and reduce to parametric {\em quasi-variational inequalities} if the sets $C(p)$ are convex; see, e.g., \cite{FP,K2,L2,LR,mo,R03,R13,Y2,Y1} and the references therein for various terminology and stability results concerning parametric systems of type \eqref{qvi}.

We introduce in Section~6 a new second-order qualification condition under the name of ``general uniform second-order sufficient condition" (GUSOSC) and show that it completely characterizes Lipschitzian full stability of solutions to \eqref{qvi} under the simultaneous validity of the partial Mangasarian-Fromovitz and constant rank constraint qualifications for \eqref{ine}. If both these constraint qualifications are replaced by the stronger linear independence constraint qualification (LICQ) for the partial gradients of the active constraints in \eqref{ine} at the reference point, then the new GUSOSC reduces to the known ``general strong second-order sufficient condition" (GSSOSC) from \cite{K2}, a slight modification for variational conditions of Robinson's strong second-order sufficient condition \cite{R} in parametric nonlinear programming with ${\cal C}^2$-smooth data. In this way we arrive at a rather surprising result that the GSSOSC characterizes the Lipschitzian full stability notion from Section~4 in the framework of parametric variational conditions in \eqref{qvi}, \eqref{ine} (and thus for parametric quasi-variational inequalities if the sets $C(p)$ are convex) under the LICQ assumption imposed in \eqref{ine}. Since Lipschitzian full stability in \eqref{qvi}, \eqref{ine} implies the local single-valuedness and Lipschitz continuity of the solution map for this system, our new GUSOSC gives a sufficient condition for the latter conventional properties. We present an example in Section~6 showing that the GUSOSC holds and thus ensures the local single-valuedness and Lipschitz continuity of the solution map to \eqref{qvi} with linear constraints in \eqref{ine} while the well-recognized ``strong coherent orientation condition" known to be sufficient for these properties in such a setting \cite{FP} fails.

The developed notions of Lipschitzian and H\"olderian full stability and their characterizations for general parametric variational systems in their specifications obtained in this paper have undoubted potentials for applications to particular variational models governed by ordinary and partial differential equations as well as to qualitative and numerical aspects of optimization, equilibria, and control. These issues and related topics will be considered in our future research.\vspace*{0.05in}

Notation and terminology of the paper are standard in variational analysis and generalized differentiation; cf. \cite{M1,rw}. Unless otherwise stated,  throughout the paper we assume that the {\em decision} space $X$ is {\em Hilbert} being identified with its dual space $X^*$. As usual, the symbol $\la\cdot,\cdot\ra$ indicates the canonical pairing in $X$ with the norm $\|x\|:=\sqrt{\la x,x\ra}$ while the symbol $\st{w}{\to}$ signifies the weak convergence in $X$. We denote by $\B$ the closed unit ball in the space in question, and thus $\B_\eta(x):=x+\eta\B$ stands for the closed ball centered at $x$ with radius $\eta>0$. Given a set-valued mapping $F\colon X\tto X$ from $X$ into itself $(X=X^*)$, the symbol
\begin{equation}\label{pk}
\begin{array}{ll}
\disp\Limsup_{x\to\ox}F(x):=\Big\{v\in X\Big|&\exists\;\mbox{ sequences }\;x_k\to\ox,\;v_k\st{w}{\to}v\;\mbox{ such that}\\
&v_k\in F(x_k)\;\mbox{ for all }\;k\in\IN:=\{1,2,\ldots\}\Big\}
\end{array}
\end{equation}
signifies the {\em sequential Painlev\'e-Kuratowski outer/upper limit} of $F(x)$ as $x\to\ox$. As stated as the beginning, the {\em parameter} space $(P,d)$ is {\em metric}, and we denote by $\B_\eta(p):=\{q\in P|\;d(q,p)\le\eta\}$ the closed ball centered at $p$ with radius $\eta>0$. The closed ball in the product space $X\times P$ is referred as $\B_\eta(x,p):=\B_\eta(x)\times\B_\eta(p)$. Recall finally that the symbols $x\st{f}{\to}\ox$ and $x\st{\O}{\to}\ox$ for a function $f\colon X\to\oR$ and a set $\O\subset X$ indicate that $x\to\ox$ with $f(x)\to f(\ox)$ and $x\in\O$, respectively.

\section{Generalized Differentiation and Preliminary Material}
\setcounter{equation}{0}

First we present here the generalized differential constructions for function, sets, and set-valued mappings widely implied in the paper; see \cite{M1,rw} for more details. Given an extended-real-valued function $f\colon X\to\oR$ on a Hilbert space $X$, supposed unless otherwise stated that it is {\em l.s.c.} around the reference points. The {\em regular subdifferential} (known also as the Fr\'echet or viscosity subdifferential) of $f$ at the point $\ox$ from $\dom f:=\{x\in X|\;f(x)<\infty\}$ is
\begin{eqnarray}\label{2.1}
\Hat\partial f(\ox):=\Big\{v\in X\Big|\;\liminf_{x\to\ox}\frac{f(x)-f(\ox)-\la v,x-\ox\ra}{\|x-\ox\|}\ge 0\Big\},
\end{eqnarray}
while the {\em limiting subdifferential} (known also as the basic subdifferential and as the Mordukhovich subdifferential) and the {\em singular subdifferential} (known also as the horizon subdifferential) of $f$ at $\ox\in\dom f$  are defined via the sequential outer limit \eqref{pk}  by
\begin{eqnarray}\label{2.2}
\partial f(\ox):=\Limsup_{x\st{f}{\to}\ox}\Hat\partial f(x)\;\mbox{ and }\;\partial^\infty f(\ox):=\Limsup_{x\st{f}{\to}\ox,\lm\dn 0}\lm\Hat\partial f(x),
\end{eqnarray}
respectively. If $f$ is convex, both regular and limiting subdifferentials reduce to the subdifferential of convex analysis. Furthermore, we have $\partial^\infty f(\ox)=\{0\}$ if $f$ is locally Lipschitzian around $\ox$.

Given a nonempty set $\O\subset X$ locally closed around $\ox\in\O$, the regular and limiting {\em normal cones} to $\O$ at $\ox\in\O$ are defined, respectively, by
\begin{eqnarray}\label{nc}
\Hat N_\O(\ox):=\Hat\partial\delta_\O(\ox)\;\mbox{ and }\;N_\O(\ox):=\partial\delta_\O(\ox)
\end{eqnarray}
via the corresponding subdifferential constructions \eqref{2.1} and \eqref{2.2} applied to the indicator function $\delta_\O(x)$ of $\O$ equal to $0$ for $x\in\O$ and to $\infty$ otherwise.

Let $F:X\tto Y$  be  a set-valued mapping between Hilbert spaces with the domain $\dom F:=\{x\in X|\;F(x)\ne\emp\}$ and the graph $\gph F:=\{(x,y)\in X\times Y|\;y\in F(x)\}$. Assume that $\gph F$ is locally closed around $(\ox,\oy)\in\gph F$ and define the {\em regular coderivative} and the {\em limiting coderivative} of $F$ at $(\ox,\oy)$ by using the corresponding normal cone \eqref{nc} to the graph of $F$ by, respectively,
\begin{eqnarray}
\Hat D^*F(\ox,\oy)(w)&:=&\big\{z\in X\big|\;(z,-w)\in\Hat N_{{\rm gph} F}(\ox,\oy)\big\},\quad w\in Y,\label{2.3}\\
 D^*F(\ox,\oy)(w)&:=&\big\{z\in X\big|\;(z,-w)\in N_{{\rm gph}F}(\ox,\oy)\big\},\quad w\in Y.\label{2.4}
\end{eqnarray}
When $F$ is single-valued around $\ox$, we skip $\oy=F(\ox)$ from the coderivative notation. It has been strongly recognized that the coderivatives (\ref{2.3}) and (\ref{2.4}) are appropriate tools for the study and characterizations of well-posedness and sensitivity in nonlinear and variational analysis; see \cite[Chapter~4]{M1} for more details and references. Recall to this end the {\em Lipschitz-like} (known also as pseudo-Lipschitz or Aubin) property of $F$ around $(\ox,\oy)\in\gph F$ defined in the case when $X$ is a metric space with metric $d$ while $Y$ is a normed space as follows: there are neighborhoods $U$ of $\ox$ and $V$ of $\oy$ as well as a constant $\ell>0$ such that
\begin{equation}\label{2.5}
F(x)\cap V\subset F(u)+\ell d(x,u)\B\;\mbox{ for all }\;x,u\in U.
\end{equation}
When both $X$ and $Y$ are finite-dimensional, the Lipschitz-like property of $F$ admits a {\em pointwise characterization} known as the coderivative/Mordukhovich criterion
\begin{eqnarray}\label{cod-cr}
D^*F(\ox,\oy)(0)=\{0\}
\end{eqnarray}
used in what follows; see \cite[Corollary~5.4]{m93} and \cite[Theorem~9.40]{rw}. We refer the reader to \cite[Theorem~4.10]{M1} for additional conditions ensuring the validity of criterion \eqref{cod-cr} in infinite dimensions.

It is worth mentioning that the limiting constructions above enjoy comprehensive pointwise calculus rules (``full calculus") while their regular counterparts obey the so-called ``fuzzy calculus" rules involving neighborhood points. Both of these calculi are based on {\em variational/extremal principles} of variational analysis; see \cite{M1,rw} and the references therein.

Let us emphasize that the main results of this paper are expressed in terms of {\em second-order} subdifferential (or generalized Hessian) constructions defined in the vein proposed by the first author (see, e.g., \cite{M0,M1}) as coderivatives of first-order subgradient mappings. Actually we use here two second-order constructions of this type generated correspondingly by the regular \eqref{2.3} and limiting \eqref{2.4} coderivative of the limiting subdifferential \eqref{2.2} for extended-real-valued functions.\vspace*{0.03in}

The next definition of single-valued localization plays an important role in this paper. Note that our definition is slightly different from the one in \cite[p.~4]{DR}, where it is not required that the single-valued localization $\Hat T$ has the full domain in $U$. The reason why the condition $\dom\Hat T=U$ is needed is that most of the time we consider single-valued and {\em continuous} localization.

\begin{Definition} {\bf (single-valued localizations).}\label{locali} Let $T:X\tto Y$ be a set-valued mapping between metric spaces, and let $(\ox,\oy)\in \gph T$. We say that $T$ admits a {\sc single-valued localization} around $(\ox,\oy)$ if there is a neighborhood $U\times V\subset X\times Y$ of $(\ox,\oy)$ such that the mapping $\Hat T:U\to V$ defined by $\gph\Hat T:=\gph T\cap(U\times V)$ is single-valued on $U$ with $\dom\Hat T=U$. In this case we say that $\Hat T$ is a single-valued localization of $T$ relative to $U\times V$. If in addition $\Hat T$ is $($Lipschitz$)$ continuous on $U$, then we say that $T$ admits a {\sc $($Lipschitz$)$ continuous single-valued localization} around $(\ox,\oy)$, i.e., $\Hat T$ is a $($Lipschitz$)$ continuous single-valued localization of $T$ relative to $U\times V$.
\end{Definition}

The mapping $\Hat T:U\tto V$ from Definition~\ref{locali} is simply called a {\em localization} of $T$ relative to $U\times V$ if its uniqueness is not postulated. Note that in the above definitions we can equivalently replace the product neighborhood  $U\times V$ by any open set $W$ around $(\ox,\oy)$. Indeed, define
$$
{\rm Pr}_X(W):=\big\{x\in X\big|\;\exists\,y\in Y\;\mbox{ with }\;(x,y)\in W\big\}
$$
and observe that this set is a neighborhood of $\ox$. Denoting now $\Hat T\colon{\rm Pr}_X(W)\tto Y$ by
$$
\Hat T(x):=\big\{y\in Y\big|\;y\in T(x),\;(x,y)\in W\big\},
$$
we get that $\gph\Hat T=\gph T\cap W$, i.e., $\hat T$ is a localization of $T$ relative to $W$.

Recall finally that $T\colon X\tto Y$ is {\em strongly metrically regular} around $(\ox,\oy)$ with modulus $\kk>0$ if the inverse mapping $T^{-1}\colon Y\tto X$ admits a Lipschitz continuous single-valued localization around $(\oy,\ox)$ with constant $\kk>0$; see, e.g., \cite{DR,FP} for further details.

\section{Characterizations of Local Strong Maximal Monotonicity}
\setcounter{equation}{0}

In this self-contained section we study some local monotonicity properties of set-valued operators in Hilbert spaces and obtain complete {\em coderivative characterizations} of {\em local strong maximal monotonicity}, which is actually behind the quantitative full stability notions for PVS \eqref{1.1} and their specifications studied in the subsequent sections of the paper.

Given a Hilbert space $X$, recall that an operator $T:X\tto X$ is (globally) {\em monotone} if
\[
\la v_1-v_2,u_1-u_2\ra\ge 0\quad\mbox{whenever}\quad(u_1,v_1),(u_2,v_2)\in\gph T.
\]
The monotone operator $T$ is {\em maximal monotone} if $\gph T=\gph S$ for any monotone operator $S\colon X\tto X$ satisfying the inclusion $\gph T\subset\gph S$. The next definition presents several types of {\em local} monotonicity (cf.\ \cite{P,PR2}) considered in this section.

\begin{Definition}{\bf (local monotonicity)}\label{mono} Let $T\colon X\tto X$, and let $(\ox,\ov)\in\gph T$. We say that:

$\bullet$ $T$ is {\sc locally monotone} around $(\ox,\ov)$ if there is a neighborhood $U\times V$ of $(\ox,\ov)$ such that
\begin{eqnarray}\label{lm}
\la v_1-v_2,u_1-u_2\ra\ge 0\quad\mbox{for all}\quad(u_1,v_1),(u_2,v_2)\in\gph T\cap(U\times V).
\end{eqnarray}
$T$ is {\sc locally maximally monotone} around $(\ox,\ov)$ if there is a neighborhood $U\times V$ of $(\ox,\ov)$ such that \eqref{lm} holds and that $\gph T\cap(U\times V)=\gph S\cap(U\times V)$ for any monotone operator $S:X\tto X$ satisfying $\gph T\cap(U\times V)\subset\gph S$.

$\bullet$  $T$ is $($locally$)$ {\sc hypomonotone} around $(\ox,\ov)$ if there exist a neighborhood $U\times V$ of this point and some positive number $r$ such that
\begin{eqnarray}\label{hm}
\la v_1-v_2,u_1-u_2\ra\ge-r\|u_1-u_2\|^2\quad\mbox{for all}\quad(u_1,v_1),(u_2,v_2)\in\gph T\cap(U\times V).
\end{eqnarray}

$\bullet$  $T$ is {\sc locally strongly monotone} around $(\ox,\ov)$ with modulus $\kk>0$  if there exists a neighborhood $U\times V$ of $(\ox,\ov)$ such that
\begin{eqnarray}\label{sm}
\la v_1-v_2,u_1-u_2\ra\ge\kk\|u_1-u_2\|^2\quad\mbox{for all}\quad(u_1,v_1),(u_2,v_2)\in\gph T\cap(U\times V).
\end{eqnarray}
Finally, $T$ is {\sc locally strongly maximally monotone} around $(\ox,\ov)$ with modulus $\kk>0$ if there exists a neighborhood $U\times V$ such that \eqref{sm} holds  and that $\gph T\cap(U\times V)=\gph S\cap(U\times V)$ for any monotone operator $S:X\tto X$ satisfying $\gph T\cap(U\times V)\subset\gph S$.
\end{Definition}

First we briefly discuss local hypomonotonicity. It is shown in \cite{P,PR1,PR2,rw} that this class of operators is rather broad including locally monotone operators and Lipschitzian mappings, limiting subgradient mappings for continuously prox-regular functions considered in Section~4, etc. The next proposition presents two useful relationships involving hypomonotonicity and localization.

\begin{Proposition}{\bf (hypomonotonicity and single-valued localization).}\label{pro1} Let $T_1,T_2:X\tto X$ be set-valued mappings with $(\ox,\ov_1)\in\gph T_1$ and $(\ox,\ov_2)\in\gph T_2$. The following assertions hold:

{\bf (i)} If $T_1$ admits a Lipschitz continuous single-valued localization around $(\ox,\ov_1)$, then $T_1$ is hypomonotone around $(\ox,\ov_1)$.

{\bf (ii)} If both $T_1$ and $T_2$ are hypomonotone around $(\ox,\ov_1)\in\gph T_1$ and $(\ox,\ov_2)\in\gph T_2$, respectively, then $T:=T_1+T_2$ is also hypomonotone around $(\ox,\ov)$ with $\ov:=\ov_1+\ov_2$ provided that $T_1$ has a continuous single-valued localization around $(\ox,\ov_1)$.
\end{Proposition}
{\bf Proof.} To verify {\bf(i)}, employ Definition~\ref{locali} and find a neighborhood $U\times V$ of $(\ox,\ov_1)$ such that the single-valued localization $S\colon U\to V$ with $\gph S=\gph T_1\cap(U\times V)$ is Lipschitz continuous on $U$ with some constant $\ell>0$. For any $(u_1,v_1),(u_2,v_2)\in\gph T_1\cap(U\times V)=\gph S$ we have
\[
\la v_1-v_2,u_1-u_2\ra\ge -\|v_1-v_2\|\cdot\|u_1-u_2\|=-\|S(u_1)-S(u_2)\|\cdot\|u_1-u_2\|\ge -\ell\|u_1-u_2\|^2,
\]
which shows by definition \eqref{hm} that $T_1$ is hypomonotone around $(\ox,\ov_1)$.

Next we justify the hypomonotonicity calculus rule formulated in {\bf(ii)}. In the notation of (ii), let $\vt$ be a continuous single-valued localization of $T_1$ relative to some neighborhood $U\times V$ of $(\ox,\ov_1)$. By the hypomonotonicity of $T_i$, $i=1,2$, around the corresponding points we find neighborhoods $U_1\times V_1\subset U\times V$ of $(\ox,\ov_1)$ and $U_2\times V_2$ of $(\ox,\ov_2)$ and two numbers $r_1,r_2>0$ such that
\begin{eqnarray}\label{hy}
\la v_1-v_2,u_1-u_2\ra\ge-r_i\|u_1-u_2\|^2\quad\mbox{for all}\quad(u_1,v_1),(u_2,v_2)\in\gph T_i\cap(U_i\times V_i)
\end{eqnarray}
for $i=1,2$. Since $\vt\colon U\to V$ is continuous, there is a neighborhood $\Tilde U\times\Tilde V$ of $(\ox,\ov)$ such that $\Tilde U\subset U_1\cap U_2$, $\vt(\Tilde U)\subset V_1$, and $\Tilde V-\vt(\Tilde U)\subset V_2$. Now pick any $(u_1,v_1),(u_2,v_2)\in\gph T\cap(\Tilde U\times\Tilde V)$ and define $v^\prime_i:=\vt(u_i)=T_1(u_i)\cap V_1$, $i=1,2$. We get $v_i-v^\prime_i\in T_2(u_i)\cap\big(\Tilde V-\vt(\Tilde U)\big)\subset T_2(u_i)\cap V_2$ for $i=1,2$. It follows from \eqref{hy} that
\begin{eqnarray*}
\la v^\prime_1-v^\prime_2,u_1-u_2\ra\ge-r_1\|u_1-u_2\|^2\quad\mbox{and}\quad \la(v_1-v^\prime_1)-(v_2-v^\prime_2),u_1-u_2\ra\ge-r_2\|u_1-u_2\|^2.
\end{eqnarray*}
Adding these two inequalities side by side gives us that
\[
\la v_1-v_2,u_1-u_2\ra\ge-(r_1+r_2)\|u_1-u_2\|^2,
\]
which ensures the claimed hypomonotonicity of $T$ around $(\ox,\ov)$.\endproof

The next result crucial in what follows describes local strong maximal monotonicity of operators in Hilbert spaces via single-valued localizations of their inverses.

\begin{Lemma}{\bf (local strong maximal monotonicity via single-valued localizations).}\label{pro2} Let $T:X\tto X$ be a set-valued mapping with $(\ox,\ov)\in\gph T$. The following are equivalent:

{\bf (i)} $T$ is locally strongly maximally monotone around $(\ox,\ov)$ with modulus $\kk>0$.

{\bf (ii)} $T$ is locally strongly monotone around $(\ox,\ov)$ with modulus $\kk$ and the inverse mapping $T^{-1}$ admits a Lipschitz continuous single-valued localization around $(\ov,\ox)$.

{\bf (iii)} The mapping $T^{-1}$ admits a single-valued localization $\vt$ relative to a neighborhood $V\times U$ of $(\ov,\ox)$  such that for all $v_1,v_2\in V$ we have the estimate
\begin{eqnarray}\label{3.5}
\big\|(v_1-v_2)-2\kk\big[\vt(v_1)-\vt(v_2)\big]\big\|\le \|v_1-v_2\|.
\end{eqnarray}
Consequently, if $T$ is locally strongly maximally monotone around $(\ox,\ov)$, then $T$ is strongly metrically regular around $(\ox,\ov)$ with modulus $\kk^{-1}$.
\end{Lemma}
{\bf Proof.} To verify implication {\bf (i)}$\Longrightarrow${\bf (ii)}, suppose that $T$ is locally strongly maximally monotone around $(\ox,\ov)$ and then find a neighborhood $U\times V$ of $(\ox,\ov)$ such that \eqref{sm} is satisfied and that the equality $\gph T\cap(U\times V)=\gph S\cap(U\times V)$ holds for any monotone operator $S:X\tto X$ satisfying $\gph T\cap(U\times V)\subset\gph S$. Denote $W:=J_\kk(U\times V)$ with $J_\kk(u,v):=(u,v-\kk u)$ for $(u,v)\in X\times X$ and note from \eqref{sm} that the set-valued mapping $F:X\tto X$ defined by $\gph F:=\gph(T-\kk I)\cap W$ is monotone. Indeed, for any $(u_i,v_i)\in\gph F$, $i=1,2$ we have $(u_i,v_i+\kk u_i)\in\gph T\cap J_\kk^{-1}(W)=\gph T\cap(U\times V)$. It follows from \eqref{sm} that
\[
\la v_1+\kk u_1-v_2-\kk u_2, u_1-u_2\ra\ge\kk\|u_1-u_2\|^2,
\]
which implies that $\la v_1-v_2,u_1-u_2\ra\ge 0$ and thus justifies the monotonicity of $F$. Accordingly, there exists a maximal monotone operator $R:X\tto X$ extending $F$ via Zorn's lemma (see, e.g., \cite[Theorem~20.21]{BC}), which means that $\gph F\subset\gph R$ and that $R$ is maximal monotone. It yields
\[
\gph(F+\kk I)\cap(U\times V)=\gph T\cap(U\times V)\subset\gph (R+\kk I).
\]
The local maximality of $T$ relative to $U\times V$ and the monotonicity of $R+\kk I$ ensure the representation $\gph T\cap(U\times V)=\gph(R+\kk I)$, and thus we have
\begin{eqnarray}\label{3.6}
\gph T^{-1}\cap(V\times U)=\gph (R+\kk I)^{-1}\cap(V\times U).
\end{eqnarray}
Applying the classical Minty theorem tells us that $\dom(R+\kk I)^{-1}=X$ and that the mapping $(R+\kk I)^{-1}$ is Lipschitz continuous on $X$. Combining this with \eqref{3.6} gives us that $V_1:=(R+\kk I)(U)\cap V$ is a neighborhood of $\ov$. Moreover, it follows from \eqref{3.6} that  $T^{-1}(v)=(R+\kk I)^{-1}(v)$ for all $v\in V_1$. Thus the localization $S:V_1\to U$ with $\gph S=\gph T^{-1}\cap(V_1\times U)$ is single-valued and Lipschitz continuous in $V_1$. This justifies implication {\bf(i)}$\Longrightarrow${\bf(ii)}.

To prove now implication {\bf(ii)}$\Longrightarrow${\bf(iii)}, suppose that $T$ is locally strongly monotone around $(\ox,\ov)$ and that $T^{-1}$ admits a Lipschitz continuous single-valued localization $\vt$ relative to some neighborhood $V\times U$ of $(\ov,\ox)$.  By shrinking $U,V$ if necessary, we get that condition \eqref{sm} is also valid on this neighborhood $U\times V$. For any $(v_1,u_1),(v_2,u_2)\in\gph\vt$, observe from \eqref{sm} that
\begin{eqnarray*}
\|v_1-v_2-2\kk(u_1-u_2)\|^2=\|v_1-v_2\|^2-4\kk\big[\la v_1-v_2,u_1-u_2\ra-\kk\|u_1-u_2\|^2\big]\le\|v_1-v_2\|^2,
\end{eqnarray*}
which therefore justifies assertion {\bf (iii)}.

It remains to verify implication {\bf(iii)}$\Longrightarrow${\bf(i)}. Pick any $(u_1,v_1),(u_2,v_2)\in\gph T\cap(U\times V)$, where $V\times U$ is the neighborhood of $(\ov,\ox)$ on which $T^{-1}$ admits a single-valued localization $\vt$ satisfying \eqref{3.5}. It tells us that $u_1=\vt(v_1)$, $u_2=\vt(v_2)$, and it follows from \eqref{3.5} that
\[
0\le\|v_1-v_2\|^2-\|v_1-v_2-2\kk(u_1-u_2)\|^2=4\kk\big[\la v_1-v_2,u_1-u_2\ra-\kk\|u_1-u_2\|^2\big].
\]
This clearly gives us the estimates
\begin{eqnarray}\label{3.7}
\|v_1-v_2\|\cdot\|u_1-u_2\|\ge\la v_1-v_2,u_1-u_2\ra\ge\kk\|u_1-u_2\|^2,
\end{eqnarray}
which not only verify \eqref{sm} but also show that $\vt$ is Lipschitz continuous in $V$ with constant $\kk^{-1}$. Then we deduce from \cite[Lemma~2.1]{MN} that $\vt$ is maximal monotone relative to $U\times V$, and so is $T$. This justifies {\bf(i)} and completes the proof of the equivalencies. The final consequence of the lemma follows from {\bf (ii)} by the definition of strong metric regularity.\endproof

Next we establish the main result in this section, which provides a characterization of local strong maximal monotonicity via the regular coderivative \eqref{2.3} for set-valued mappings in Hilbert spaces. The first result in this direction turns back to \cite[Theorem~2.1]{PR2}, where Poliquin and Rockafellar obtained a necessary condition for the global maximal monotonicity in terms of the limiting coderivative \eqref{2.4} in finite dimensions. More recently \cite{CT}, Chieu and Trang established necessary and sufficient coderivative conditions for global monotonicity and strong monotonicity
for single-valued and continuous mappings. Our result below gives a coderivative characterization of {\em local strong maximal} monotonicity of general {\em set-valued} mappings.

\begin{Theorem}{\bf (regular coderivative characterization of local strong maximal monotonicity of set-valued mappings).}\label{char1} Let $T:X\tto X$ be a set-valued mapping for which the set $\gph T$ is locally closed around the point $(\ox,\ov)\in\gph T$. The following are equivalent:

{\bf (i)} $T$ is locally strongly maximally monotone around $(\ox,\ov)$ with modulus $\kk>0$.

{\bf (ii)} $T$ is hypomonotone around $(\ox,\ov)$ and there exists $\eta>0$ such that
\begin{eqnarray}\label{3.9}
\la z,w\ra\ge\kk\|w\|^2\quad\mbox{for all}\quad z\in\Hat D^*T(u,v)(w),\;(u,v)\in\gph T\cap\B_\eta(\ox,\ov).
\end{eqnarray}
The conditions in {\bf (ii)} ensure the strong metric regularity of $T$ around $(\ox,\ov)$ with modulus $\kk^{-1}$.
\end{Theorem}
{\bf Proof.} We start with justifying {\bf(i)}$\Longrightarrow${\bf(ii)}. It is obvious that $T$ is hypomonotone around $(\ox,\ov)$ when {\bf (i)} is satisfied. By Lemma~\ref{pro2} there is a single-valued localization $\vt$ of $T^{-1}$ relative to some neighborhood $V\times U$ of $(\ov,\ox)$ such that inequality \eqref{3.5} holds.  Observe from \eqref{3.7}, which is a consequence of \eqref{3.5} by the proof of Lemma~\ref{pro2}, that $\vt$ is Lipschitz continuous on $V$ with modulus $\kk^{-1}$. Fix $\eta>0$ satisfying $\B_\eta(\ox,\ov)\subset U\times V$ and pick any $(u,v)\in\gph T\cap\B_\eta(\ox,\ov)$ and $z\in\Hat D^*T(u,v)(w)$. With $\ve>0$ we find, by definition \eqref{2.3} of the regular coderivative, some number $\delta>0$ such that $\B_\delta(u,v)\subset U\times V$ and that
\begin{eqnarray}\label{3.10}
\ve(\|x-u\|+\|y-v\|)\ge\la z,x-u\ra-\la w,y-v\ra\quad\mbox{for all}\quad (x,y)\in\gph T\cap\B_\delta(u,v).
\end{eqnarray}
When $t>0$ is sufficiently small, define $u_t:=\vt(v_t)$ with $v_t:=v+t(z-2\kk w)\in V$ and get from the local Lipschitz continuity of $\vt$ that $(u_t,v_t)\to(u,v)$ as $t\dn 0$. Without loss of generality, suppose that $(u_t,v_t)\in\B_\delta(u,v)$ for all $t>0$. Replacing $(x,y)$ in \eqref{3.10} by $(u_t,v_t)$ and using \eqref{3.5} yield
\begin{eqnarray}\label{3.11}
\begin{array}{ll}
\ve\big(\|u_t-u\|+\|v_t-v\|\big)&\disp\ge\la z,u_t-u\ra-\la w,v_t-v\ra\\
&\disp=\la t^{-1}(v_t-v)+2\kk w,u_t-u\ra-t\la w,z-2\kk w\ra\\
&\disp=t^{-1}\la v_t-v,u_t-u\ra+2\kk\la w,u_t-u\ra-t\la w,z\ra+2t\kk\|w\|^2\\
&\disp\ge\kk t^{-1}\|u_t-u\|^2+2\kk\la w,u_t-u\ra-t\la w,z\ra+2t\kk\|w\|^2\\
&\disp\ge\kk t^{-1}\|u_t-u\|^2-2\kk\|w\|\cdot\|u_t-u\|+t\kk\|w\|^2-t\la w,z\ra+t\kk\|w\|^2\\
&\disp\ge-t\la z,w\ra+t\kk\|w\|^2.
\end{array}
\end{eqnarray}
Since $\vt$ is Lipschitz continuous on $V$ with modulus $\kk^{-1}$, we have
\[
\begin{array}{ll}
\ve\big(\|u_t-u\|+\|v_t-v\|\big)&\disp=\ve\big(\|\vt(v_t)-\vt(v)\|+\|v_t-v\|\big)\le\ve\big(\kk^{-1}\|v_t-v\|+\|v_t-v\|\big)\\
&\disp=\ve(\kk^{-1}+1)\|v_t-v\|=\ve(\kk^{-1}+1)t\|z-2\kk w\|,
\end{array}
\]
which together with \eqref{3.11} yields $\la z,w\ra+\ve(\kk^{-1}+1)\|z-2\kk w\|\ge\kk \|w\|^2$, and so $\la z,w\ra\ge\kk\|w\|^2$ by taking $\ve\dn 0$. This ensures \eqref{3.9} and thus completes the proof of {\bf(i)}$\Longrightarrow${\bf(ii)}.\vspace*{0.05in}

To verify the converse implication {\bf(ii)$\Longrightarrow$(i)}, observe that by Lemma~\ref{pro2} we only need to show that $T^{-1}$ admits a Lipschitz continuous single-valued localization $\vt$ around $(\ov,\ox)$, which satisfies \eqref{3.5}. It is done below in the following two claims.\\[-1.5ex]

{\bf Claim~1.} {\em $T^{-1}$ admits a Lipschitz continuous localization $\vt$ around $(\ov,\ox)$}.\\[-1.5ex]

By choosing $\eta>0$ to be sufficiently small, we may always assume that the set $\gph T\cap\B_\eta(\ox,\ov)$ is closed and there is a number $r>0$ such that
\begin{eqnarray}\label{3.12}
\la v_1-v_2,x_1-x_2\ra\ge -r\|x_1-x_2\|^2\quad\mbox{for all}\quad(x_1,v_1),(x_2,v_2)\in\gph T\cap\B_\eta(\ox,\ov).
\end{eqnarray}
Pick any $s>r$ and define $J_s(u,v):=(v+su,u)$ for $(u,v)\in X\times X$. Denote further $W_s:=J_s(\B_\eta(\ox,\ov))$ and observe that ${\rm int}\,W_s=J_s({\rm int}\,\B_\eta(\ox,\ov))$ is a neighborhood of $(\ov+s\ox,\ox)$. It follows from \eqref{3.12} that for all $(v_1,x_1),(v_2,x_2)\in\gph(T+sI)^{-1}\cap W_s:=\gph F$ we have
\begin{eqnarray}\label{3.13}
\|v_1-v_2\|\cdot\|x_1-x_2\|\ge\la v_1-v_2,x_1-x_2\ra\ge (s-r)\|x_1-x_2\|^2,
\end{eqnarray}
which implies that the localization $F$ of $(T+sI)^{-1}$ is single-valued. Taking any $(v,u)\in\gph F\cap({\rm int}\,W_s)$ and $w\in\Hat D^*F(v,u)(z)$, we get that $w\in \Hat D^*(T+sI)^{-1}(v,u)(z)$ and that $-z\in\Hat D^* (T+sI)(u,v)(-w)$. It follows from the coderivative sum rule \cite[Theorem~1.62]{M1} that $-z+sw\in\Hat D^*T(u,v-su)(-w)$. Since $(u,v-s u)=J_s^{-1}(u,v)\in J^{-1}_s({\rm int}\,W_s)={\rm int}\,\B_\eta(\ox,\ov)$, we deduce from \eqref{3.9} that $\la-z+ s w,-w\ra\ge\kk\|w\|^2$, and thus
\begin{eqnarray}\label{3.14}
\|z\|\cdot\|w\|\ge\la z,w\ra\ge(\kk+s)\|w\|^2.
\end{eqnarray}
To proceed further, for any $z\in\B$ define the extended-real-valued function
\begin{eqnarray}\label{3.15}
f_z(v):=\left\{\begin{array}{ll}
\la z, F(v)\ra\quad&\mbox{if}\quad v\in\dom F,\\
\infty&\mbox{otherwise}.
\end{array}\right.
\end{eqnarray}
Since $\gph T\cap \B_\eta(\ox,\ov)$ is closed in $X\times X$, it is easy to see that $\gph F$ is also closed on $X\times X$. Let us show that $f_z$ is l.s.c.\ on $X$. Arguing by contradiction, suppose that there exist $\ve>0$ and a sequence $v_k$ converging to some $v\in X$ such that $f_z(v_k)<f(v)-\ve$. If $f_z(v)=\infty$, then we have $v\notin\dom F$ while $v_k\in\dom F$. It follows from \eqref{3.13} that $\|F(v_k)-F(v_j)\|\le(s-r)^{-1}\|v_k-v_j\|$, and so $F(v_k)$ is a Cauchy sequence. Hence the sequence $(v_k,F(v_k))\in\gph F$ converges to $(v,x)\in  \gph F$ due to the closedness of $\gph F$, which implies that $F(v)=x$ and contradicts $v\notin\dom F$. If $f_z(v)<\infty$, then \eqref{3.13} ensures that $\|F(v_k)-F(v)\|\le(s-r)\|v_k-v\|\to 0$, which contradicts $f_z(v_k)\to f_z(v)$.

Fix now a positive number $\delta<\frac{\eta}{3}$ and pick any $(u_i,v_i)\in\gph T\cap\B_{\delta}(\ox,\ov)$, $i=1,2$. Then we have $(y_i,u_i)\in\gph F$ with $y_i:=v_i+su_i$. Choosing $\ve\in(0,\delta)$ and applying the mean value inequality \cite[Corollary~3.50]{M1} to the l.s.c.\ function $f_z$ give us that
\begin{eqnarray}\label{3.16}
|f_z(y_1)-f_z(y_2)|\le\|y_1-y_2\|\sup\big\{\|w\|\;\big|\;w\in\hat\partial\la z,F\ra(y),\;y\in[y_1,y_2]+\ve\B\big\},
\end{eqnarray}
where $[y_1,y_2]:=\{ty_1+(1-t)y_2|\;t\in[0,1]\}$. For any $y\in\dom F\cap\big([y_1,y_2]+\ve\B\big)$ there are some $t\in[0,1]$ and $y_0\in\ve\B$ such that $y=ty_1+(1-t)y_2+y_0$. It follows from \eqref{3.13} that
\begin{eqnarray*}\begin{array}{ll}
\|y-\ov-s\ox\|&\disp=\|ty_1+(1-t)y_2+y_0-\ov-s\ox\|\\
\disp&=\|t(y_1-\ov-s\ox)+(1-t)(y_2-\ov-s\ox)+y_0\|\\
\disp&=\|t(v_1+su_1-\ov-s\ox)+(1-t)(v_2+su_2-\ov-s\ox)+y_0\|\\
\disp&\le t\big(\|v_1-\ov\|+s\|u_1-\ox\|\big)+(1-t)\big(\|v_2-\ov\|+s\|u_2-\ox\|\big)+\|y_0\|\\
\disp&\le t(\delta+s\delta)+(1-t)(\delta+s\delta)+\ve=(1+s)\delta+\ve<(2+s)\delta.
\end{array}
\end{eqnarray*}
We easily get from the latter estimate and \eqref{3.13} that
\begin{eqnarray}\label{3.17}
\|F(y)-\ox\|=\|F(y)-F(\ov+s\ox)\|\le (s-r)^{-1}\|y-\ov-s\ox\|\le (s-r)^{-1}(2+s)\delta.
\end{eqnarray}
Furthermore, it follows from the above that
\begin{eqnarray}\label{3.18}\begin{array}{ll}
\|y-s F(y)-\ov\|&\disp=\|y-\ov-s\ox-s(F(y)-\ox)\|\le \|y-\ox-s\ox\|+s\|F(y)-\ov\|\\
&\disp\le (2+s)\delta+s(s-r)^{-1}(2+s)\delta.
\end{array}
\end{eqnarray}
By choosing $\delta$ sufficiently small, we get from \eqref{3.17} and \eqref{3.18} that $J_s^{-1}(y, F(y))=(F(y),y-sF(y))\in \inte\B_\eta(\ox,\ov)$, which yields $(y,F(y))\in J_s\big(\inte\B_\eta(\ox,\ov)\big)=\inte W_s$. Moreover, note that
\[
\hat\partial\la z,F\ra(y)\subset\hat D^*F(y)(z)=\hat D^*F\big(y,F(y)\big)(z).
\]
Since $(y,F(y))\in\gph(T+sI)^{-1}\cap\inte W_s$, we have $\hat D^*F(y,F(y))(z)= \Hat D^*(T+sI)^{-1}(y,F(y))(z)$. This together with \eqref{3.16}, \eqref{3.15}, and \eqref{3.14} shows that
\[
|\la z,F(y_1)-F(y_2)\ra|=|f_z(y_1)-f_z(y_2)|\le\|y_1-y_2\|\cdot(\kk+s)^{-1}\|z\|\quad\mbox{for all}\quad z\in\B.
\]
It allows us to conclude that
\begin{eqnarray*}
\|u_1-u_2\|=\|F(y_1)-F(y_2)\|\le(\kk+s)^{-1}\|y_1-y_2\|=(\kk+s)^{-1}\|v_1+su_1-v_2-su_2\|,
\end{eqnarray*}
which implies in turn the inequality
\[
(\kk+s)\|u_1-u_2\|\le\|(v_1-v_2)+s(u_1-u_2)\|\le \|v_1-v_2\|+s\|u_1-u_2\|.
\]
Thus we arrive at the estimate
\begin{eqnarray}\label{3.19}
\kk\|u_1-u_2\|\le\|v_1-v_2\|\quad  \mbox{for all}\quad (u_1,v_1),(u_2,v_2)\in\gph T\cap\B_{\delta}(\ox,\ov).
\end{eqnarray}

It remains to check that the inverse mapping $T^{-1}$ admits a Lipschitzian localization around $(\ov,\ox)$. Observe to this end from \eqref{3.9} that
\begin{eqnarray*}
\|z\|\ge\kk\|w\|\quad\mbox{for all}\quad z\in\Hat D^*T(u,v)(w),\;(u,v)\in\gph T\cap\B_\eta(\ox,\ov).
\end{eqnarray*}
It follows from the regular coderivative criterion in \cite[Theorem~4.7]{M1} that $T^{-1}$ is Lipschitz-like around $(\ov,\ox)$ with some modulus $\ell>0$. By definition \eqref{2.5} we find $\nu>0$ such that the inclusion $\B_{\ell\nu}(\ox)\times\B_\nu(\ov)\subset\B_\delta(\ox,\ov)$ holds and that
\[
\ox\in T^{-1}(v)+\ell\|v-\ov\|\quad\mbox{for all}\quad v\in \B_\nu(\ov),
\]
which implies in turn that $T^{-1}(v)\cap\inte\B_{\ell\nu}(\ox)\ne\emp$ for all $v\in\inte\B_\nu(\ov)$. Defining finally the mapping $\vt$ from $\inte\B_\nu(\ov)$ to $\inte\B_{\ell\nu}(\ox)$ by $\gph\vt:=\gph T^{-1}\cap(\inte\B_\nu(\ov)\times\inte\B_{\ell\nu}(\ox))$, we have $\dom\vt=\inte\B_\nu(\ov)$, and it follows from \eqref{3.19} that $\vt$ is single-valued. This together with \eqref{3.19} shows that $\vt$ is locally Lipschitz continuous with modulus $\kk$. \\[-1.5ex]

{\bf Claim~2.} {\em The single-valued localization $\vt$ of $T^{-1}$ defined in Claim~1 satisfies inequality \eqref{3.5}}.\\[-1.5ex]

For any $z\in\B$ we define $\xi_z(v):=\la z,v-2\kk\vt(v)\ra$, $v\in\B_\nu(\ov)$. Fix $\al,\beta>0$ with $\al+\beta<\nu$ and $v_1,v_2\in \B_\al(\ov)$. Similarly to \eqref{3.16} we get from the mean value inequality \cite[Corollary~3.50]{M1} that
\begin{equation}\label{3.20}
|\xi_z(v_1)-\xi_z(v_2)|\le \|v_1-v_2\|\sup\big\{\|w\|\;\big|\; w\in\hat\partial\xi_z(v),\;v\in[v_1,v_2]+\beta\B\big\}.
\end{equation}
Since $v\in\inte\B_\nu(\ov)$ for each $v\in[v_1,v_2]+\beta\B$, it tells us that
\begin{equation*}
w\in\hat\partial\xi_z(v)\subset z-2\kk\hat D^*\vt(v)(z)=z-2\kk\hat D^* T^{-1}(v)(z),
\end{equation*}
which yields $(2\kk)^{-1}(z-w)\in\hat D^* T^{-1}(v)(z)$, or equivalently $-z\in\hat D^*T(\vt(v),v)((2\kk)^{-1}(w-z))$. Then we deduce from \eqref{3.9} the inequality
\[
\la-z,(2\kk)^{-1}(w-z)\ra\ge k\|(2\kk)^{-1}(w-z)\|^2,
\]
which implies in turn that $\|w\|\le\|z\|$. This together with \eqref{3.20} ensures that
\[
|\xi_z(v_1)-\xi_z(v_2)|\le\|v_1-v_2\|\cdot\|z\|\quad \mbox{for all}\quad z\in\B.
\]
Remembering the definition of $\xi_z$, we arrive at the estimate
\[
\|v_1-v_2-2\kk(\vt(v_1)-\vt(v_2))\|\le\|v_1-v_2\|\quad\mbox{whenever}\quad v_1,v_2\in\B_\al(\ov),
\]
which verifies \eqref{3.5} and thus justifies Claim~2. This completes the proof of the theorem by combining the results given in Claim~1 and Claim~2.
\endproof

Note that the aforemention strong metric regularity of $T$ can be characterized by using the strict graphical derivative in finite-dimensions; see \cite[Definition~9.53 and Theorem~9.54]{rw}. Our result above provides a verifiable sufficient condition \eqref{3.9} for this property in terms of the regular coderivative under the hypomonotonicity assumption. However, the main trust of Theorem~\ref{char1} is a characterizations of the local strong maximal monotonicity property, which significantly supersedes strong metric regularity and is needed in what follows.\vspace*{0.03in}

Next we derive from Theorem~\ref{char1} a {\em pointwise} characterization of the local strong maximal monotonicity property for single-valued Lipschitzian mappings in finite-dimensional spaces via the limiting coderivative \eqref{2.4}. This is actually a natural extension of the classical result stated that a $\mathcal{C}^1$-smooth mapping $F:\R^n\to\R^n$ is locally strongly monotone around $\ox$ provided that $\nabla F(\ox)$ is positive-definite. Note further that in the latter case the local maximality and hypomonotonicity of $F$ are automatic due to the Lipschitz continuity of this mapping.

\begin{Corollary}{\bf (limiting coderivative characterization of local strong monotonicity for Lipschitz continuous mappings).}\label{coro1} Let $X$ be a finite-dimensional space, and let $T:X\to X$ be a single-valued mapping Lipschitz continuous around $\ox\in\dom T$. The following are equivalent:

{\bf (i)} $T$ is locally strongly monotone around $(\ox,T(\ox))$.

{\bf (ii)} $D^*T(\ox)$ is positive-definite in the sense that
\begin{eqnarray}\label{psd}
\la z,w\ra>0\quad\mbox{whenever}\quad z\in D^*T(\ox)(w),\;w\ne 0.
\end{eqnarray}
\end{Corollary}
{\bf Proof.} It suffices to check that condition \eqref{psd} is equivalent to \eqref{3.9} under the assumptions made. By passing to the limit it is easy to derive implication \eqref{3.9}$\Longrightarrow$\eqref{psd}. To justify the converse implication, we argue by contradiction and suppose that \eqref{psd} is satisfied while \eqref{3.9} is not. This gives us a sequence $(u_k,z_k,w_k)$ such that $u_k\to\ox$, $z_k\in\Hat D^*T(u_k)(w_k)$, and $\la z_k,w_k\ra<\frac{1}{k}\|w_k\|^2$, which implies that $w_k\ne 0$. Define further $\ow_k:=\frac{w_k}{\|w_k\|}$ and $\oz_k:=\frac{z_k}{\|w_k\|}$. Since $T$ is Lipschitz continuous around $\ox$ with some modulus $\ell$, we have $\|\oz_k\|\le\ell\|\ow_k\|=\ell$ for sufficiently large $k$. By passing to subsequences, assume without loss of generality that $\ow_k\to\ow$ with $\|\ow\|=1$ and $\oz_k\to\oz$ as $k\to\infty$. It follows from definition \eqref{2.4} of the limiting coderivative that $\oz\in D^*T(\ox)(\ow)$. Furthermore, by $\la\oz_k,\ow_k\ra<\frac{1}{k}$ we get the inequality $\la\oz,\ow\ra\le 0$, which contradicts the positive-definiteness condition \eqref{psd} and thus completes the proof of the corollary.\endproof

As a direct consequence of Corollary~\ref{coro1}, observe that condition \eqref{psd} is sufficient for the strong metric regularity of $T$ around $(\ox,T(\ox))$. It has been proved by Kummer \cite{Ku} that the latter property can be characterized by using Thibault's strict derivative \cite{th}. To this end we emphasize again that our positive-definiteness coderivative criterion \eqref{psd} characterizes essentially more specific property of local strong maximal monotonicity of our main interest here.\vspace*{0.03in}

Finally in this section, we formulate a conjecture for which the affirmative answer is achieved in the cases presented in Corollary~\ref{coro1} as well as in Corollary~\ref{coro4} given in the next section.

\begin{Conjecture}{\bf (limiting coderivative characterization of local strong maximal monotonicity for set-valued mappings).}\label{con} Let $X$ be a finite-dimensional space, and let $T:X\tto X$ be a set-valued mapping with closed graph around $(\ox,\ov)\in\gph T$. The following are equivalent:

{\bf (i)} $T$ is locally strongly maximally monotone around $(\ox,\ov)$.

{\bf (ii)} $T$ is hypomonotone around $(\ox,\ov)$ and $D^*T(\ox,\ov)$ is positive-definite in the sense that
\begin{eqnarray}\label{pd}
\la z,w\ra>0\quad \mbox{whenever}\quad z\in D^*T(\ox,\ov)(w),\;w\ne 0.
\end{eqnarray}
\end{Conjecture}

Note that implication (i)$\Longrightarrow$(ii) of this conjecture follows from Theorem~\ref{char1} by using the limiting procedure. However, the converse implication would be more interesting. In Corollary~\ref{coro4} the reader can find the justification of this conjecture in the important set-valued case including subgradient mappings generated by a major and fairly broad class of extended-real-valued functions.

\section{Characterizations of Full Stability in Variational Systems}
\setcounter{equation}{0}

In this section we turn to the main subject of the paper concerning full quantitative stability of the parametric variational systems (PVS) given by
\begin{eqnarray}\label{VS}
 v\in f(x,p)+\partial_x g(x,p)\;\mbox{ for }\;x\in X,\;p\in P
\end{eqnarray}
with the Hilbert decision space $X$ and the metric parameter space $(P,d)$, where $f:X\times P\to X$, $g:X\times P\to\oR$, and $\partial_x g$ stands for the partial limiting subdifferential of the function $g$ with respect to the variable $x$. Denote $g_p(\cdot):=g(\cdot,p)$ and observe that $\partial_x g(x,p)=\partial g_p(x)$ for all $(x,p)\in X\times P$.

Fix $\ov\in f(\ox,\op)+\partial_x g(\ox,\op)$ and consider the solution map $S\colon X\times P\tto X$ to \eqref{VS} defined by
\begin{eqnarray}\label{ss}
S(v,p):=\big\{x\in X\big|\;v\in f(x,p)+\partial_x g(x,p)\big\}\;\mbox{ with }\;\ox\in S(\ov,\op).
\end{eqnarray}
The underlying goal of this section is to introduce and efficiently characterize the following new notions of H\"olderian and Lipschitzian full stability for PVS \eqref{VS}.

\begin{Definition}{\bf (H\"olderian and Lipschitzian full stability of parametric variational systems).}\label{fs} Given $\ox\in S(\ov,\op)$ from \eqref{ss}, we say that:

{\bf (i)} $\ox$ is a {\sc H\"olderian fully stable} solution to PVS \eqref{VS} corresponding to the parameter pair $(\ov,\op)$ if the solution map \eqref{ss} admits a single-valued localization $\vt$ relative to some neighborhood $V\times Q\times U$ of $(\ov,\op,\ox)$ such that for any $(v_1,p_1), (v_2,p_2)\in V\times Q$ we have
\begin{eqnarray}\label{4.7}
\big\|(v_1-v_2)-2\kk[\vt(v_1,p_1)-\vt(v_2,p_2)]\big\|\le\|v_1-v_2\|+\ell d(p_1,p_2)^{\frac{1}{2}}
\end{eqnarray}
with some positive constants $\kk$ and $\ell$.

{\bf (ii)} $\ox$ is a {\sc Lipschitzian fully stable} solution to PVS \eqref{VS} corresponding to the parameter pair $(\ov,\op)$ if the solution map \eqref{ss} admits a single-valued localization $\vt$ relative to some neighborhood $V\times Q\times U$ of $(\ov,\op,\ox)$ such that for any $(v_1,p_1), (v_2,p_2)\in V\times Q$ we have
\begin{eqnarray}\label{mp1}
\big\|(v_1-v_2)-2\kk\big[\vt(v_1,p_1)-\vt(v_2,p_2)\big]\big\|\le\|v_1-v_2\|+\ell d(p_1,p_2)
\end{eqnarray}
with some positive constants $\kk$ and $\ell$.
\end{Definition}

It is easy to see that the stability notions defined above imply the {\em local single-valuedness} and {\em Lipschitz} (resp.\ {\em H\"older}) {\em continuity} of the solution map $S$, which are conventional definitions of quantitative stability of perturbed variational systems discussed, e.g., in \cite{L2,LR,Ro,R03,Y2} for the Lipschitzian case and in \cite{Y1} for the H\"olderian one. However, full stability from Definition~\ref{fs} is {\em much stronger} than these conventional notions even in very simple settings. Consider, e.g., $g=0$ and $f\colon\R^2\to\R^2$ in \eqref{VS} given by $f(x):=(x_1,-x_2)$. It is obvious that $f^{-1}$ is single-valued and Lipschitz continuous around $(0,0)$ while the Lipschitzian full stability property \eqref{mp1} fails.

More generally, we can observe to this end that when the basic parameter $p$ is omitted in \eqref{VS}, both stability conditions \eqref{4.7} and \eqref{mp1} reduce to the one in \eqref{3.5}, which is equivalent to the {\em local strong maximal monotonicity} of the mapping $T:=f+\partial g$ by Lemma~\ref{pro2}. Thus for $g=0$ and $f=f(x)$ the full stability conditions of Definition~\ref{fs} amount to saying that $f$ is strongly monotone around $\ox$, which essentially supersedes the Lipschitz continuity of its inverse. We will see in what follows that the methods developed in Section~3 to characterize this monotonicity concept, particularly Theorem~\ref{char1} and its proof, play a crucial role in deriving efficient second-order criteria for the full stability notions from Definition~\ref{fs}.

Furthermore, it will be shown below that the introduced notions of full stability for PVS \eqref{VS} are {\em equivalent} in the case of $f=0$ to the  corresponding full stability definitions for {\em local minimizers} associated with $g\colon X\times P\to\oR$, which are initiated in \cite{LPR} for the Lipschitzian version and then extended in \cite{MN} to H\"olderian full stability in optimization. Let us recall these definitions.

We say that $\ox\in X$ is a {\em Lipschitzian fully stable local minimizer} associated with $g\colon X\times P\to\oR$ relative to $\op\in P$ with $(\ox,\op)\in\dom g$ and some ``tilt" parameter $\ov\in X$ if there exist positive numbers $\kappa,\ell,\gg$ and a neighborhood $V\times Q$ of $(\ov,\op)$ such that the argminimum mapping
\[
(v,p)\mapsto M_\gg(v,p):={\rm argmin}\,\big\{g(x,p)-\la v,x\ra\big|\;x\in\B_\gg(\ox)\big\}
\]
is single-valued on $V\times Q$ with $M_\gg(\ov,\op)=\ox$ and satisfies the Lipschitz condition
\begin{eqnarray}\label{LS}
\|M_\gg(v_1,p_1)-M_\gg(v_2,p_2)\|\le\kk\|v_1-v_2\|+\ell d(p_1,p_2)\;\mbox{ for all }\;v_1,v_2\in V,\;p_1,p_2\in Q.
\end{eqnarray}

The point $\ox$ is a {\em H\"olderian fully stable local minimizer} associated to of $g$ relative to $\op$ and $\ov$ if there exist positive numbers $\kappa,\ell,\gg$ such that the argminimum mapping $M_\gg$ is single-valued on some neighborhood $V\times Q$ of $(\ov,\op)$ with $M_\gg(\ov,\op)=\ox$ and satisfies the condition
\begin{eqnarray}\label{HFS}
\|M_\gg(v_1,p_1)-M_\gg(v_2,p_2)\|\le\kk\|v_1-v_2\|+\ell d(p_1,p_2)^{\frac{1}{2}}\;\;\mbox{for all}\;\;v_1,v_2\in V,\;p_1,p_2\in Q.
\end{eqnarray}
These notions have studied intensively in both Lipschitzian \cite{LPR,mos,mrs,ms} and H\"olderian \cite{MN,MNR} cases for various constrained optimization problems in finite and infinite dimensions with deriving second-order characterizations of fully stable local minimizers therein.

We can see that definitions \eqref{LS} and \eqref{HFS} of full stability for local minimizers are formulated differently in comparison with our new Definition~\ref{fs} of full stability for parametric variational systems; the former ones essentially exploit specific features of scalar optimization. The equivalence between these types of full stability in the optimization framework is a nontrivial fact that follows from the criteria of full stability for parametric variational systems obtained below and those established earlier for local minimizers. In this way, the full stability conditions in Definition~\ref{fs} can be treated as new full stability characterizations for local minimizers in parametric optimization. On the other hand, the derivation of the full stability criteria for parametric variational systems given below includes the construction (by using the coderivative conditions for local strong maximal monotonicity from Section~3 and advanced techniques of nonsmooth variational analysis) of auxiliary problems of extended-real-valued optimization and employing second-order characterizations of fully stable local minimizers; see the proof of Theorem~\ref{Holder}.\vspace{0.03in}

To proceed with the formulation and proof of our main result in this section, we first specify the class of functions $g$ from \eqref{VS} used in our analysis. In fact, it is the major and fairly large collection of extended-real-valued functions employed in second-order variational analysis and parametric optimization. In the parametric framework this class of functions has been defined in \cite{LPR} by extending the corresponding nonparametric versions widely used in variational analysis, optimization, and their applications in both finite and infinite dimensions; see, e.g., \cite{BT1,PR1,rw} and the references therein. Given $g:X\times P\to\oR$ finite at $(\ox,\op)$ with $\hat v:=\ov-f(\ox,\op)\in\partial_x g(\ox,\op)$, we say following \cite{LPR} that $g$ is {\em prox-regular} in $x$ at $\ox$ for $\hat v$ with {\em compatible parameterization} by $p$ at $\op$ if there are neighborhoods $U$ of $\ox$, $V$ of $\hat v$, and $Q$ of $\op$ along with numbers $\ve>0$ and $r>0$ such that
\begin{eqnarray}\label{par-prox}
\begin{array}{ll}
\quad\quad g(x,p)\ge g(u,p)+\la v,x-u\ra-\frac{r}{2}\|x-u\|^2\;\mbox{ for all }\;x\in U,\\
\quad\mbox{when }\;v\in\partial_x g(u,p)\cap V,\;u\in U,\;p\in Q,\;\mbox{ and }\;g(u,p)\le g(\ox,\op)+\ve.
\end{array}
\end{eqnarray}
Further, $g$ is {\em subdifferentially continuous} in $x$ at $\ox$ for $\hat v$ with compatible parameterization by $p$ at $\op$ if the mapping $(x,p,v)\mapsto f(x,p)$ is continuous relative to the subdifferential graph $\gph\partial_x g$ at $(\ox,\op,\hat v)$. For simplicity we call $g$ to be {\em parametrically continuously prox-regular} at $(\ox,\op)$ for $\hat v$ when $g$ is simultaneously prox-regular and subdifferentially continuous at $\ox$ for $\hat v$ with compatible parameterization by $p$ at $\op$. In this case inequality ``$g(u,p)\le g(\ox,\op)+\ve$'' can be omitted in \eqref{par-prox}.\vspace*{0.03in}

Throughout this section we impose the following {\em standing assumptions} on the data of \eqref{VS}:\vspace*{0.03in}

{\bf (A1)} $f$ is differentiable with respect to $x$ around $(\ox,\op)$ uniformly in $x$ and the partial Jacobian $\nabla_x f$ is continuous at $(\ox,\op)$. Furthermore, $f$ is Lipschitz continuous with respect to $p$ uniformly in $x$ around $(\ox,\op)$, i.e., there exist a neighborhood $U\times Q$ of $(\ox,\op)$ and a constant $L>0$ such that
\begin{equation}\label{4.4b}
\|f(x,p_1)-f(x,p_2)\|\le Ld(p_1,p_2)\quad\mbox{for all}\quad x\in U,\,p_1,p_2\in Q.
\end{equation}

{\bf (A2)} $g$ is  parametrically continuously prox-regular at $(\ox,\op)$ for $\hat v$.

{\bf (A3)} The following {\em basic constraint qualification} (BCQ) holds at $(\ox,\op)$:
\begin{eqnarray}\label{bcq}
\mbox{the mapping}\;\; p\mapsto\epi g(\cdot,p)\;\;\mbox{is Lipschitz-like around}\;\big(\op,(\ox,g(\ox,\op))\big).
\end{eqnarray}

Note that assumption {\bf(A1)} for $f$ is classical in the study of generalized equations and turns back to the landmark paper by Robinson \cite{Ro}. It easily follows from the mean value theorem that such a mapping $f$ is also Lipschitz continuous around $(\ox,\op)$, i.e., with no change the notation in comparison with \eqref{4.4b}, there exist a neighborhood $U\times Q$ of $(\ox,\op)$ and a constant $L>0$ for which
\begin{equation}\label{4.6}
\|f(x_1,p_1)-f(x_2,p_2)\|\le L\big[\|x_1-x_2\|+d(p_1,p_2)\big]\;\mbox{ whenever }\;(x_1,p_1),(x_2,p_2)\in U\times Q.
\end{equation}
Observe also that in the case when both spaces $X$ and $P$ are finite-dimensional, BCQ from {\bf (A3)} can be equivalently described by the implication \begin{eqnarray}\label{bcq2}
(0,q)\in\partial^\infty g(\ox,\op)\Longrightarrow q=0,
\end{eqnarray}
which follows from the coderivative criterion \eqref{cod-cr} for the Lipschitz-like property of the epigraphical mapping in \eqref{bcq}; see \cite{LPR}. It is shown in Sections~5 and 6 that both assumptions {\bf(A2)} and {\bf(A3)} naturally hold for important special classes of parametric variational systems in finite and infinite dimensions. It is worth also mentioning furthermore that when the parameter $p$ is not present, we have assumption {\bf (A3)} to fulfill automatically, {\bf (A1)} means that $f$ is continuously differentiable around $\ox$, and {\bf (A2)} reduces to the continuous prox-regularity of $g$ at $\ox$ for $\hat v$.

Now we are ready to formulate the main results of this section giving a second-order characterization of H\"older full stability for general parametric variational systems in Hilbert spaces.\vspace*{0.03in}

\begin{Theorem}{\bf (second-order characterization of H\"olderian full stability for PVS).}\label{Holder} Let $\ox\in S(\ov,\op)$ from \eqref{ss} under our standing assumptions. Consider the following two statements:

{\bf (i)} $\ox$ is a H\"olderian fully stable solution of PVS \eqref{VS} corresponding to the parameter pair $(\ov,\op)$ with the moduli $\kk,\ell>0$ taken from \eqref{4.7}.

{\bf (ii)} There exist numbers $\eta,\kk_0>0$ such that whenever $(u,p,v)\in\gph\partial_x g\cap\B_\eta(\ox,\op,\ov)$ we have
\begin{eqnarray}\label{4.8}
\la\nabla_x f(\ox,\op)w,w\ra+\la z,w\ra\ge\kk_0\|w\|^2\quad\mbox{for all}\quad z\in(\Hat D^*\partial g_p)(u,v)(w),\;w\in X.
\end{eqnarray}
Then {\bf (i)} implies {\bf (ii)} with constant $\kk_0$ that can be chosen smaller than but arbitrarily close to $\kk$. Conversely, the validity of {\bf (ii)} ensures that {\bf (i)} holds, where $\kk$ can be chosen smaller but arbitrarily close to $\kk_0$. Consequently, \eqref{4.8} implies that the solution map \eqref{ss} admits a single-valued and H\"older continuous localization $\vt$ relative to a neighborhood $V\times Q\times U$ of $(\ov,\op,\ox)$, i.e., for any $(v_1,p_1),(v_2,p_2)\in V\times Q$ we have
\begin{eqnarray}\label{4.9}
\big\|\vt(v_1,p_1)-\vt(v_2,p_2)\big\|\le\frac{1}{\kk}\|v_1-v_2\|+\frac{\ell}{2\kk}d(p_1,p_2)^{\frac{1}{2}}.
\end{eqnarray}
\end{Theorem}

As significant steps in the proof of Theorem~\ref{Holder}, we begin with deriving the following two lemmas, which are of their own interest while being largely employed in the subsequent parts of this section. The first lemma establishes a certain ``time propagation" of the aforementioned full stability properties in the case of  linearized PVS of type \eqref{VS}.

\begin{Lemma}{\bf (propagation of full stability for linearized PVS).}\label{lm1} Denote $A:=\nabla_x f(\ox,\op)$ and consider the one-parametric family of operators
$$
A_t:=\disp\frac{1}{2}\Big(A+A^*\Big)+tB\;\mbox{ with }\;B:=A-A^*\;\mbox{ and }\;t\ge 0,
$$
where $A^*$ is the adjoint operator of $A$. Define further the set-valued mapping $G_t:X\times P\tto X$ by
\begin{eqnarray}\label{gt}
G_t(v,p):=\big\{x\in X\big|\;v\in f(\ox,\op)+A_t(x-\ox)+\partial_x g(x,p)\big\}\quad\mbox{for all}\quad(v,p)\in X\times P.
\end{eqnarray}
The following two assertions are satisfied:

{\bf (i)} Suppose that $G_\tau$ for some $\tau\ge 0$ has a single-valued localization $\vt_\tau$ relative to a neighborhood $V\times Q\times U$ of $(\ov,\op,\ox)$ such that for any $ (v_1,p_1),(v_2,p_2)\in V\times Q$ we have
\begin{eqnarray}\label{Hol1}
\big\|(v_1-v_2)-2\kk\big[\vt_\tau(v_1,p_1)-\vt_\tau(v_2,p_2)\big]\big\|\le\|v_1-v_2\|+\ell d(p_1,p_2)^\frac{1}{2}
\end{eqnarray}
with some $\ell>0$. Then $G_t$ also admits a single-valued localization $\vt_t$ relative to a neighborhood $V_1\times Q_1\times U_1\subset V\times Q\times U$ of $(\ov,\op,\ox)$ such that for any $(v_1,p_1),(v_2,p_2)\in V_1\times Q_1$ we have
\begin{eqnarray}\label{Hol2}
\big\|(v_1-v_2)-2\kk\big[\vt_t(v_1,p_1)-\vt_t(v_2,p_2)\big]\big\|\le\|v_1-v_2\|+2\ell d(p_1,p_2)^\frac{1}{2}
\end{eqnarray}
whenever $t\in\big[\tau,\tau+\frac{\kk}{2\|B\|}\big)$ under the convention that $1/0:=\infty$.

{\bf (ii)} Suppose that $G_\tau$ for some $\tau\ge 0$ has a single-valued localization $\vt_\tau$ relative to a neighborhood $V\times Q\times U$ of $(\ov,\op,\ox)$ such that for any $ (v_1,p_1),(v_2,p_2)\in V\times Q$ we have
\begin{eqnarray}\label{Lip1}
\big\|(v_1-v_2)-2\kk\big[\vt_\tau(v_1,p_1)-\vt_\tau(v_2,p_2)\big]\big\|\le\|v_1-v_2\|+\ell d(p_1,p_2)
\end{eqnarray}
with some $\ell>0$. Then $G_t$ also admits a single-valued localization $\vt_t$ relative to a neighborhood  $V_1\times Q_1\times U_1\subset V\times Q\times U$ of $(\ov,\op,\ox)$ satisfying the condition
\begin{eqnarray}\label{Lip2}
\big\|(v_1-v_2)-2\kk\big[\vt_t(v_1,p_1)-\vt_t(v_2,p_2)\big]\big\|\le\|v_1-v_2\|+2\ell d(p_1,p_2)
\end{eqnarray}
 for any $(v_1,p_1),(v_2,p_2)\in V_1\times Q_1$, provided that $t\in\big[\tau,\tau+\frac{\kk}{2\|B\|}\big)$ under the convention above.
\end{Lemma}
{\bf Proof.} To justify assertion {\bf (i)}, observe first from \eqref{Hol1} that
\[
2\kk\|\vt_\tau(v_1,p_1)-\vt_\tau(v_2,p_2)\|-\|v_1-v_2\|\le\|v_1-v_2\|+\ell d(p_1,p_2)^\frac{1}{2}
\]
for all $(v_1,p_1),(v_2,p_2)\in V\times Q$. This clearly implies the estimate
\begin{eqnarray}\label{Holder1}
\|\vt_\tau(v_1,p_1)-\vt_\tau(v_2,p_2)\|\le\ell_1\|v_1-v_2\|+\ell_2 d(p_1,p_2)^\frac{1}{2}\quad \mbox{with}\quad\ell_1:=\frac{1}{\kk},\; \ell_2:=\frac{\ell}{2\kk}.
\end{eqnarray}
Fixing $t\in\big[\tau,\tau+\frac{\kk}{2\|B\|}\big)$ from the formulation of the lemma, we obviously get $\ve:=r(1-\ell_1(t-\tau)\|B\|)\in(0,r]$ for any $r>0$. Choose now $r>0$ so small that $\B_{r\ell_1}(\ox)\subset U$, $\B_r(\ov)\subset V$, and $\B_s(\op)\subset Q$ with $s:={\frac{\ve^2\ell^2_1}{4\ell_2^2}}$. Define the mapping $\vt_t$ by
\[
\gph\vt_t:=\gph G_t\cap\big(V_1\times Q_1\times U_1)
\]
with $V_1:={\rm int}\,\B_\frac{\ve}{2}(\ov)\subset\B_r(\ov)\subset V$, $Q_1:={\rm int}\,\B_s(\op)\subset Q$, $U_1={\rm int}\,\B_{r\ell_1}(\ox)\subset U$ and show that it is single-valued. To proceed, for arbitrary parameters $v\in V_1$ and $p\in Q_1$ we form the single-valued mapping $T:\B_{r\ell_1}(\ox)\to X$ by
\begin{eqnarray}\label{T}
T(x):=\vt_\tau\big(v-(t-\tau)B(x-\ox),p\big)\quad\mbox{for all}\quad x\in\B_{r\ell_1}(\ox)
\end{eqnarray}
and claim that $T$ has the full domain $\dom T=\B_{r\ell_1}(\ox)$. Indeed, it is easy to see that
\[
\|v-(t-\tau)B(x-\ox)-\ov\|\le \|v-\ov\|+(t-\tau)\|B\|\cdot\|x-\ox\|<\ve+(t-\tau)r\ell_1\|B\|=r\;\mbox{ for }\;x\in \B_{r\ell_1}(\ox),
\]
which amounts to $v-(t-\tau)B(x-\ox)\in\B_r(\ov)\subset V$. We get from the definition of $\vt_\tau$ that the values $T(x)=\vt_\tau(v-(t-\tau)B(x-\ox),p)$ are well defined, and thus $\dom T=\B_{r\ell_1}(\ox)$. Observe further that $u\in\vt_t(v,p)$ if and only if $u=T(u)$, i.e., $u$ is a fixed point of \eqref{T}.  Moreover, the H\"older continuity of $\vt_\tau$ in \eqref{Holder1} yields the estimates
\begin{equation}\label{4.17}
\|T(\ox)-\ox\|=\|\vt_\tau(v,p)-\vt_\tau(\ov,\op)\|\le \ell_1\|v-\ov\|+\ell_2d(p,\op)^\frac{1}{2}<\frac{\ell_1\ve}{2}+\ell_2s^{\frac{1}{2}}=\frac{\ell_1\ve}{2}+\frac{\ell_1\ve}{2}= \ell_1\ve.
\end{equation}
It follows from \eqref{Holder1} that for any $x_1,x_2\in\B_{r\ell_1}(\ox)$ we have
\begin{eqnarray}\begin{array}{ll}\label{4.18}
\|T(x_1)-T(x_2)\|&\disp=\|\vt_\tau(v-(t-\tau)B(x_1-\ox),p)-\vt_\tau(v-(t-\tau)B(x_2-\ox),p)\|\\
&\disp\le\ell_1(t-\tau)\|B(x_1-x_2)\|\le\ell_1(t-\tau)\|B\|\cdot\|x_1-x_2\|.
\end{array}
\end{eqnarray}
 Pick now any $x\in\B_{r\ell_1}(\ox)$ and deduce from \eqref{4.17} and \eqref{4.18} that
\begin{eqnarray}\label{4.19}
\begin{array}{ll}
\|T(x)-\ox\|&\disp\le\|T(x)-T(\ox)\|+\|T(\ox)-\ox\|< \ell_1(t-\tau)\|B\|\cdot\|x-\ox\|+\ell_1\ve\\
&\disp\le\ell_1(t-\tau)\|B\|r\ell_1+\ell_1r(1-\ell_1(t-\tau)\|B\|)=r\ell_1.
\end{array}
\end{eqnarray}
Since $\ell_1(t-\tau)\|B\|<1$ due to the choice of $t$, this allows us to apply the classical contraction principle to the mapping $T$ satisfying \eqref{4.18} and \eqref{4.19} and thus to find a unique point $u\in\B_{r\ell_1}(\ox)$ with $T(u)=u$. It follows from \eqref{4.19} that $u\in U_1$. Thus $\vt_t(v,p)$ is a {\em singleton} whenever $(v,p)\in V_1\times Q_1$.\\[-1.5ex]

To complete the proof of {\bf(i)}, it remains to show that $\vt_t$ satisfies \eqref{Hol2}. Indeed, picking any $(v_i,p_i)\in V_1\times Q_1$ and denoting $u_i:=\vt_t(v_i,p_i)$, $i=1,2$, we have $u_i=\vt_\tau(v_i-(t-\tau)B(u_i-\ox),p_i)$. Define also $u_3:=\vt_t(v_1,p_2)$, which means that $u_3=\vt_\tau(v_1-(t-\tau)B(u_3-\ox),p_2)$. Employing  \eqref{Hol1} for two pairs $(v_1-(t-\tau)B(u_3-\ox),p_2)$ and $(v_2-(t-\tau)B(u_2-\ox),p_2)$ yields
\begin{eqnarray*}
\|v_1-v_2-(t-\tau)B(u_3-u_2)-2\kk(u_3-u_2)\|\le\|v_1-v_2-(t-\tau)B(u_3-u_2)\|,
\end{eqnarray*}
which ensures in turn the conditions
\[\begin{array}{ll}
0&\disp\ge 4\kk^2\|u_3-u_2\|^2-4\kk\la v_1-v_2-(t-\tau)B(u_3-u_2),u_3-u_2\ra\\
&\disp= 4\kk^2\|u_3-u_2\|^2-4\kk\la v_1-v_2,u_3-u_2\ra+4\kk(t-\tau)\la B(u_3-u_2),u_3-u_2\ra\\
&\disp=4\kk^2\|u_3-u_2\|^2-4\kk\la v_1-v_2, u_3-u_2\ra\\
&\disp=\|v_1-v_2-2\kk(u_3-u_2)\|^2-\|v_1-v_2\|^2,
\end{array}
\]
where the second equality is valid since for any $x\in X$ we have
\[
\la Bx,x\ra=\la Ax,x\ra-\la A^*x,x\ra=\la x,A^*x\ra-\la A^*x,x\ra=0.
\]
This gives us the norm relationship
\begin{eqnarray}\label{4.20}
\|v_1-v_2\|\ge\|v_1-v_2-2\kk(u_3-u_2)\|.
\end{eqnarray}
Then using \eqref{Hol1} for two pairs $(v_1-(t-\tau)B(u_1-\ox),p_1)$ and $(v_1-(t-\tau)B(u_3-\ox),p_2)$ yields
\begin{eqnarray*}
\|-(t-\tau)B(u_1-u_3)-2\kk(u_1-u_3)\|\le\|(t-\tau)B(u_1-u_3)\|+\ell d(p_1,p_2)^\frac{1}{2}
\end{eqnarray*}
with noting again that $\vt_\tau(v_1-(t-\tau)B(u_1-\ox),p_1)=\vt(v_1,p_1)=u_1$ and $\vt_\tau(v_1-(t-\tau)B(u_3-\ox),p_2)=\vt_t(v_1,p_2)=u_3$.
Thus we arrive at the lower distance estimate
\[
\ell d(p_1,p_2)^\frac{1}{2}\ge 2(\kk-(t-\tau)\|B\|)\|u_1-u_3\|\ge \kk \|u_1-u_3\|,
\]
where the last inequality is due to the choice of $t$. This tells us together with \eqref{4.20} that
\[\begin{array}{ll}
\|v_1-v_2\|&\disp\ge\|v_1-v_2-2\kk(u_1-u_2)\|-2\kk\|u_1-u_3\|\\
&\disp\ge\|v_1-v_2-2\kk(u_1-u_2)\|-2\ell d(p_1,p_2)^\frac{1}{2}.
\end{array}
\]
 It verifies \eqref{Hol2} and completes the proof of {\bf (i)}.\vspace*{0.03in}

The proof of {\bf (ii)} is quite similar. The only differences needed therein are the change of the neighborhood $Q_1$  above and the replacement of $d(p_1,p_2)^\frac{1}{2}$ by $d(p_1,p_2)$. Now we choose $Q_1:={\rm int}\,\B_s(\op)\subset Q$ with $s:=\frac{\ve\ell_1}{2\ell_2}$ for $r>0$ sufficiently small.  This allows us to show that the mapping $T$ in \eqref{T} also satisfies \eqref{4.17}, \eqref{4.18}, and \eqref{4.19} for any $v\in V_1$ and $p$ belonging to the new neighborhood $Q_1$. The rest of the proof follows the lines in the proof of {\bf(i)}.\endproof

The next lemma shows how to pass, after the parameter propagation of Lemma~\ref{lm1}, from single-valued localizations of the linearized variational systems $G_{1/2}$ from \eqref{gt} satisfying \eqref{Hol1} and \eqref{Lip1} to the corresponding single-valued localizations of the solution map $S$ to the original PVS \eqref{VS}. The linearization results of this type go back to Robinson \cite{Ro} in the case of local Lipschitz continuity for generalized equations with parameter-independent set-valued parts.

\begin{Lemma}{\bf (single-valued localizations of solutions maps to nonlinear PVS).}\label{lm2} In the setting of Lemma~{\rm\ref{lm1}} the following assertions hold:

{\bf (i)} Consider the mapping $G_\tau$ from \eqref{gt} with $\tau=\frac{1}{2}$ and suppose that it admits a single-valued localization $\vt_\tau$ relative to a neighborhood $V\times Q\times U$ of $(\ov,\op,\ox)$ satisfying \eqref{Hol1} with some moduli $\kk,\ell>0$. Then for any $\ve\in(0,\kk)$ the solution map $S$ from \eqref{ss} also admits a single-valued localization $\vt$ relative to some neighborhood $V_1\times Q_1\times U_1\subset V\times Q\times U$ of $(\ov,\op,\ox)$ such that whenever $(v_1,p_1),(v_2,p_2)\in V_1\times Q_1$ we have the estimate
\begin{eqnarray}\label{4.22a}
\big\|v_1-v_2-2(\kk-\ve)\big[\vt(v_1,p_1)-\vt(v_2,p_2)\big]\big\|\le\|v_1-v_2\|+(\ell+2L\sqrt{2\ve})d(p_1,p_2)^\frac{1}{2},
\end{eqnarray}
where the Lipschitz constant $L>0$ is taken from \eqref{4.4b}.

{\bf (ii)} Suppose that in the setting of {\bf(i)} the mapping $G_\tau$ with $\tau=\frac{1}{2}$ has a single-valued localization $\vt_\tau$ relative to a neighborhood $V\times Q\times U$ of $(\ov,\op,\ox)$ satisfying \eqref{Lip1} with moduli $\kk,\ell>0$. Then for any $\ve\in(0,\kk)$ the solution map $S$ also admits a single-valued localization $\vt$ relative to some neighborhood $V_1\times Q_1\times U_1\subset V\times Q\times U$ of $(\ov,\op,\ox)$ such that we have
\begin{eqnarray}\label{4.23a}
\big\|v_1-v_2-2(\kk-\ve)\big[\vt(v_1,p_1)-\vt(v_2,p_2)\big]\big\|\le\|v_1-v_2\|+(\ell+2L)d(p_1,p_2)
\end{eqnarray}
for all pairs $(v_1,p_1)$ and $(v_2,p_2)$ from $V_1\times Q_1$.
\end{Lemma}
{\bf Proof.} To justify {\bf (i)}, assume that $G_\tau$ with $\tau=\frac{1}{2}$ in \eqref{gt} admits a single-valued localization $\vt_\tau$ relative to a neighborhood $V\times Q\times U$ of $(\ov,\op,\ox)$ satisfying \eqref{Hol1} with some moduli $\kk,\ell>0$. Thus we also have the Holder continuity of $\vt_\tau$ in \eqref{Holder1} with $\ell_1:=\frac{1}{\kk}$ and $\ell_2:=\frac{\ell}{2\kk}$. Furthermore, observe from the construction of $G_\tau$ in \eqref{gt} and the structure of $A_\tau$ that
\[
G_{\tau}(v,p)=\big\{x\in X \big|\;v\in f(\ox,\op)+A(x-\ox)+\partial_x g(x,p)\big\}\quad\mbox{whenever}\quad (v,p)\in X\times P.
\]
 Define $r(x,p):=f(\ox,\op)+A(x-\ox)-f(x,p)$ for all $(x,p)\in X\times P$ and note that $x\in S(v,p)$ if and only if $x\in G_\tau(v+r(x,p),p)$. It follows from the standing assumption {\bf (A1)} that for any $\ve\in(0,\kk)=(0,\ell_1^{-1})$ there exist $\rho,\eta\in(0,\ve)$ with $\B_\rho(\ox)\times\B_\eta(\ov)\times\B_\eta(\op)\subset U\times V\times Q$ such that $v+r(x,p)\in U$ for $(x,v,p)\in \B_\rho(\ox)\times\B_\eta(\ov)\times \B_\eta(\op)$ and that
\begin{eqnarray}\label{4.23}\begin{array}{ll}
&\disp\|f(x_1,p)-f(x_2,p)-A(x_1-x_2)\|\le \ve\|x_1-x_2\|\quad\mbox{and}\\
&\disp\ell_1(\| f(\ox,\op)-f(\ox,p)\|+\|v-\ov\|)+\ell_2 d(p,\op)^\frac{1}{2}\le (1-\ell_1\ve)\rho
\end{array}
\end{eqnarray}
whenever $(x_1,v,p),(x_2,v,p)\in\B_\rho(\ox)\times\B_\eta(\ov)\times\B_\eta(\op)$. Fix further $(v,p)\in\B_\eta(\ov)\times\B_\eta(\op)$ and observe that the mapping $\Phi(x):=\vt_\tau(v+r(x,p),p)$ is well defined on $\B_\rho(\ox)$. Pick any vectors $x_1, x_2\in\B_\rho(\ox)$ and deduce from \eqref{Holder1} and \eqref{4.23} that
\begin{eqnarray}\begin{array}{ll}\label{4.24}
\|\Phi(x_1)-\Phi(x_2)\|&\disp\le \ell_1\|r(x_1,p)-r(x_2,p)\|\\
&\disp= \ell_1\|f(x_1,p)-f(x_2,p)-A(x_1-x_2)\|\le \ell_1\ve\|x_1-x_2\|
\end{array}
\end{eqnarray}
ensuring by $\ell_1\ve<1$ the contraction condition for $\Phi$. It also follows from \eqref{Holder1} and \eqref{4.23} that
\[\begin{array}{ll}
\|\Phi(\ox)-\ox\|&\disp=\|\vt_\tau(v+r(\ox,p),p)-\vt_\tau(\ov,\op)\|\le\ell_1\big(\|r(\ox,p)\|+\|v-\ov\|\big)+\ell_2\disp d(p,\op)^\frac{1}{2}\\
&\disp=\ell_1\big(\|f(\ox,\op)-f(\ox,p)\|+\|v-\ov\|\big)+\ell_2 d(p,\op)^\frac{1}{2}\le(1-\ell_1\ve)\rho,
\end{array}
\]
which together with the last estimate in \eqref{4.24} implies the relationships
\[
\|\Phi(x)-\ox\|\le\|\Phi(x)-\Phi(\ox)\|+\|\Phi(\ox)-\ox\|\le\ell_1\ve\|x-\ox\|+(1-\ell_1\ve)\rho\le\rho\;\;\mbox{for all}\;\;x\in \B_\rho(\ox).
\]
Combining it with \eqref{4.24} tells us that there is a unique fixed point $x$ of $\Phi$ due to the classical contraction principle. Denote this fixed point by $x(v,p):=\vt_\tau(v+r(x(v,p),p),p)$ for each $v,p\in\B_\eta(\ov)\times\B_\eta(\op)$ and then define a localization $\vt$ of $S$ by
$$
\gph\vt:=\gph S\cap\big(\inte\B_\eta(\ov)\times\inte\B_\eta(\op)\times\inte\B_\rho(\ox)\big),
$$
which is single-valued with $\vt(v,p)=x(v,p)$, $(v,p)\in\inte\B_\eta(\ov)\times\inte\B_\eta(\op)$. Picking further any $(v_i,p_i)\in\inte\B_\eta(\ov)\times\inte\B_\eta(\op),\;i=1,2$, and denoting $u_i:=\vt(v_i,p_i)$ and $u_3:=\vt(v_1,p_2)$, we get
$$
u_i=\vt_\tau\big(v_i+r(u_i,p_i),p_i\big),\;i=1,2,\;\mbox{ and }\;u_3=\vt_\tau\big(v_1+r(u_3,p_2),p_2\big).
$$
Applying \eqref{Hol1} to the pairs $(v_1+r(u_3,p_2),p_2)$ and $(v_2+r(u_2,p_2),p_2)$ gives us the estimate
\[
\|v_1-v_2+[r(u_3,p_2)-r(u_2,p_2)]-2\kk(u_3-u_2)\|\le \|v_1-v_2+[r(u_3,p_2)-r(u_2,p_2)]\|,
\]
which implies in turn the following relationships:
\[
\begin{array}{ll}
0&\disp\le\|v_1-v_2+[r(u_3,p_2)-r(u_2,p_2)]\|^2-\|v_1-v_2+(r(u_3,p_2)-r(u_2,p_2))-2\kk(u_3-u_2)\|^2\\
&\disp=4\kk\la v_1-v_2+[r(u_3,p_2)-r(u_2,p_2)],u_3-u_2\ra-4\kk^2\|u_3-u_2\|^2\\
&\disp=4\kk\la v_1-v_2,u_3-u_2\ra-4\kk^2\|u_3-u_2\|^2-4\kk\la f(u_3,p_2)-f(u_2,p_2)-A(u_3-u_2),u_3-u_2\ra\\
&\disp\le 4\kk\la v_1-v_2,u_3-u_2\ra-4\kk^2\|u_3-u_2\|^2+4\kk\|f(u_3,p_2)-f(u_2,p_2)-A(u_3-u_2)\|\cdot\|u_3-u_2\|\\
&\disp\le 4\kk\la v_1-v_2,u_3-u_2\ra-4\kk^2\|u_3-u_2\|^2+4\kk\ve\|u_3-u_2\|^2\\
&\disp=4\kk\la v_1-v_2,u_3-u_2\ra-4\kk(\kk-\ve)\|u_3-u_2\|^2,
\end{array}
\]
where the last inequality follows from \eqref{4.23}. This tells us that
\[
0\le\la v_1-v_2,u_3-u_2\ra-(\kk-\ve)\|u_3-u_2\|^2,
\]
and therefore we arrive at the estimate
\begin{eqnarray}\label{4.26}
\|v_1-v_2-2(\kk-\ve)(u_3-u_2)\|\le\|v_1-v_2\|.
\end{eqnarray}
Further, employing \eqref{Hol1} to the pairs $(v_1+r(u_3,p_2),p_2)$ and $(v_1+r(u_1,p_1),p_1)$ gives us
\[\begin{array}{ll}
\|r(u_3,p_2)-r(u_1,p_1)\|+\ell d(p_1,p_2)^\frac{1}{2}&\disp\ge\|[r(u_3,p_2)-r(u_1,p_1)]-2\kk(u_3-u_1)\|\\
&\disp\ge 2\kk\|u_3-u_1)\|-\|r(u_3,p_2)-r(u_1,p_1)\|,
\end{array}
\]
which readily implies the inequalities
\[\begin{array}{ll}
\ell d(p_1,p_2)^\frac{1}{2}&\disp\ge 2\kk\|u_3-u_1\|-2\|r(u_3,p_2)-r(u_1,p_1)\|\\
&\disp\ge 2\kk\|u_3-u_1\|-2\|r(u_3,p_2)-r(u_1,p_2)\|-2\|r(u_1,p_2)-r(u_1,p_1)\|\\
&\disp\ge 2\kk\|u_3-u_1\|-2\|f(u_1,p_2)-f(u_3,p_2)-A(u_1-u_3)\|-2\|f(u_1,p_1)-f(u_1,p_2)\|\\
&\disp\ge 2\kk\|u_3-u_1\|-2\ve\|u_3-u_1\|-2Ld(p_1,p_2),
\end{array}\]
where the last one is a consequence of \eqref{4.23} and \eqref{4.4b}. Thus we get
\[
2(\kk-\ve)\|u_3-u_1\|\le\ell d(p_1,p_2)^\frac{1}{2}+2Ld(p_1,p_2)\le\ell d(p_1,p_2)^\frac{1}{2}+2L\sqrt{2\eta}d(p_1,p_2)^\frac{1}{2}.
\]
This together with \eqref{4.26} ensures that
\[\begin{array}{ll}
\|v_1-v_2-2(\kk-\ve)(u_1-u_2)\|&\disp\le\|v_1-v_2-2(\kk-\ve)(u_3-u_2)\|+2(\kk-\ve)\|u_3-u_1\|\\
&\disp\le\|v_1-v_2\|+(\ell+2L\sqrt{2\eta})d(p_1,p_2)^\frac{1}{2}\\
&\disp\le\|v_1-v_2\|+(\ell+2L\sqrt{2\ve})d(p_1,p_2)^\frac{1}{2},
\end{array}
\]
which verifies \eqref{4.22a} and completes the proof of assertion {\bf(i)}. The proof of {\bf (ii)} is quite similar with replacing $d(p_1,p_2)^\frac{1}{2}$ by $d(p_1,p_2)$ in the arguments above.\endproof

Now we are ready to proceed with the proof of the main Theorem~\ref{Holder} by using the obtained lemmas, coderivative criteria of local strong maximal monotonicity from Section~3, and characterizations of H\"olderian full stability of local minimizers established in \cite{MN}.\vspace*{0.05in}

{\bf Proof of Theorem \ref{Holder}.} First suppose that {\bf (i)} is satisfied, i.e., the H\"olderian full stability of $\ox\in S(\ov,\op)$ from Definition~\ref{fs}(i) holds. Fix $p\in Q$ and get from \eqref{4.7} with $\vt_p(v):=\vt(v,p)$ that
\[
\big\|(v_1-v_2)-2\kk\big[\vt_p(v_1)-\vt_p(v_2)\big]\big\|\le\|v_1-v_2\|\quad\mbox{for all}\quad v_1,v_2\in V,
\]
which easily implies the {\em local strong monotonicity} condition
\begin{eqnarray}\label{4.27}
\la v_1-v_2,\vt_p(v_1)-\vt_p(v_2\ra\ge\kk\|\vt_p(v_1)-\vt_p(v_2)\|^2.
\end{eqnarray}
Define similarly $T_p(x):=f_p(x)+\partial g_p(x)$ for $x\in X$ and note that $\gph\vt_p=\gph T_p^{-1}\cap(V\times U)$. Since we automatically have the maximality in \eqref{4.27}, it follows from the coderivative characterization of Theorem~\ref{char1} that there is $\nu>0$ independent of $p$ (see the proof) such that
\begin{eqnarray}\label{4.28}
\la z,w\ra\ge\kk\|w\|^2\quad\mbox{for all}\quad z\in\Hat D^*T_p(u,v)(w),\;(u,v)\in\gph T_p\cap\B_\nu(\ox,\ov),\;w\in X.
\end{eqnarray}
By the coderivative sum rule from \cite[Theorem~1.62]{M1} we have
\[
\Hat D^*T_p(u,v)(w)=\nabla_x f(u,p)^*w+\big(\Hat D^*\partial g_p\big)\big(u,v-f(u,p)\big)(w),
\]
which says that $z\in\Hat D^*T_p(u,v)(w)$ if and only if $z=\nabla_x f(u,p)^*w+z_1$ for some $z_1\in\big(\Hat D^*\partial g_p\big)(u,v_1)(w)$ with $v_1:=v-f(u,p)$. It follows therefore that
\begin{eqnarray}\label{4.29}
\begin{array}{ll}
\la z,w\ra&\disp=\la\nabla_x f(u,p)^*w,w\ra+\la z_1,w\ra=\la\nabla_x f(u,p)w,w\ra+\la z_1,w\ra\\
&\disp\le\la\nabla_x f(\ox,\op)w,w\ra+\gg(\nu)\|w\|^2+\la z_1,w\ra
\end{array}
\end{eqnarray}
for all $p\in\B_\nu(\op)$, where the inequality holds with some $\gg(\nu)\dn 0$ as $\nu\dn 0$ due to assumption {\bf (A1)}. Hence it is possible to choose $\eta>0$ to be so small that
$$
v_1+f(u,p)\in\B_\nu(\ov)\;\mbox{ whenever }\;(u,p,v_1)\in\gph\partial_x g\cap\B_\eta(\ox,\op,\hat v).
$$
This together with \eqref{4.29} and \eqref{4.28} tells us that for any $(u,p,v_1)\in\gph\partial_x g\cap\B_\eta(\ox,\op,\hat v)$ we have
\begin{eqnarray*}
\la\nabla_x f(\ox,\op)w,w\ra+\la z_1,w\ra\ge\big(\kk-\gg(\nu)\big)\|w\|^2\;\mbox{ whenever }\;z_1\in\big(\Hat D^*\partial g_p\big)(u,v_1)(w),\;w\in X.
\end{eqnarray*}
This yields the coderivative condition \eqref{4.8} when $\nu$ is small enough to ensure that $\kk_0:=\kk-\gg(\nu)>0$ is smaller than but arbitrarily close to $\kk$.\vspace*{0.03in}

Next we verify that condition \eqref{4.8} is sufficient for the H\"olderian full stability of $\ox\in S(\ov,\op)$ while supposing that {\bf(ii)} is satisfied. Define the auxiliary l.s.c.\ function $h\colon X\times P\to\oR$ by
\begin{eqnarray}\label{hh}
h(x,p):=\la f(\ox,\op),x-\ox\ra+\frac{1}{2}\la A(x-\ox),x-\ox\ra+g(x,p)\quad\mbox{for all}\quad (x,p)\in X\times P
\end{eqnarray}
with $A=\nabla_x f(\ox,\op)$. Note that $\partial_x h(x,p)=f(\ox,\op)+\frac{1}{2}(A+A^*)(x-\ox)+\partial_x g(x,p)$ and that $\ov\in\partial_x h(\ox,\op)=f(\ox,\op)+\partial_x g(\ox,\op)$. It follows from {\bf (A2)} that $h$ is parametrically continuously prox-regular at $(\ox,\op)$ for $\ov$. For any $(u,p,v)\in \gph\partial_x h$ we get from \cite[Theorem~1.62]{M1} that
\begin{eqnarray}\label{cod1}
\big(\hat D^*\partial h_p\big)(u,v)(w)=\frac{1}{2}\Big(A+A^*\Big)w+\big(\hat D^*\partial g_p\big)\Big(u,v-f(\ox,\op)-\frac{1}{2}\Big(A+A^*\Big)(u-\ox)\Big)(w).
\end{eqnarray}
Select $\delta>0$ to be sufficiently small to ensure that
$$
\Big(u,p,v-f(\ox,\op)-\frac{1}{2}\Big(A+A^*\Big)(u-\ox)\Big)\in\B_\eta(\ox,\op,\hat v)\;\mbox{ for all }\;(u,p,v)\in\B_\delta(\ox,\op,\ov).
$$
It readily follows from \eqref{cod1} and \eqref{4.8} that
\begin{eqnarray*}
\la z,w\ra\ge\frac{1}{2}\big\la(A+A^*)w,w\big\ra+\kk_0\|w\|^2-\la A w,w\ra=\kk_0\|w\|^2\;\mbox{ if }\;z\in\big(\hat D^*\partial h_p\big)(u,v)(w)
\end{eqnarray*}
whenever $(u,p,v)\in\gph\partial_x h\cap\B_\delta(\ox,\op,\ov)$ and $w\in X$. The latter ensures by \cite[Theorem~4.7]{MN} that $\ox$ is a H\"olderian fully stable local minimizer associated with the function $h$ from \eqref{hh} relative to the parameter pair $(\ov,\op)$ as defined above. Furthermore, \cite[Theorem~3.3]{MN} applied to this function tells us that the mapping $G_0$ in \eqref{gt} admits a single-valued localization $\vt_0$ relative to some neighborhood $V_0\times Q_0\times U_0$ of $(\ov,\op,\ox)$ such that for any triple $(u,p,v)\in\gph\vt_0$ we have
\begin{eqnarray}\label{4.33}
h(x,p)\ge h(u,p)+\la v,x-u\ra+\frac{\kk_0}{2}\|x-u\|^2\quad\mbox{whenever}\quad x\in U
\end{eqnarray}
and there exists a positive number $\ell_0$ for which
\begin{eqnarray}\label{4.34}
\|\vt_0(v_1,p_1)-\vt_0(v_2,p_2)\|\le\kk^{-1}_0\|v_1-v_2\|+\ell_0 d(p_1,p_2)^\frac{1}{2}\;\mbox{ if }\;(v_1,p_1),(v_2,p_2)\in V_0\times Q_0.
\end{eqnarray}
The second-order growth condition \eqref{4.33} yields the two inequalities
\begin{eqnarray*}\begin{array}{ll}
&\disp h(u_2,p)\ge h(u_1,p)+\la v_1,u_2-u_1\ra+\frac{\kk_0}{2}\|u_2-u_1\|^2,\\\\
&\disp h(u_1,p)\ge h(u_2,p)+\la v_2,u_1-u_2\ra+\frac{\kk_0}{2}\|u_1-u_2\|^2
\end{array}
\end{eqnarray*}
for any $(u_1,p,v_1),(u_2,p,v_2)\in\gph\partial_x h\cap(U_0\times Q_0\times V_0)$.
Summing them up gives us
\begin{eqnarray}\label{m2}
\la v_1-v_2,u_1-u_2\ra\ge\kk_0\|u_1-u_2\|^2.
\end{eqnarray}
which ensures in turn the estimate
\begin{eqnarray}\label{4.35}
\|v_1-v_2-2\kk_0(\vt_0(v_1,p)-\vt_0(v_2,p))\|\le\|v_1-v_2\|\quad\mbox{for all}\quad v_1,v_2\in V,\,p\in Q.
\end{eqnarray}
Pick further any $(v_i,p_i)\in (V_0\times Q_0)$, $i=1,2$, and define $u_i:=\vt_0(v_i,p_i)$, $u_3:=\vt_0(v_1,p_2)$. Then we deduce from \eqref{4.35} and \eqref{4.34} the relationships
\[\begin{array}{ll}
\|v_1-v_2-2\kk_0(u_1-u_2)\|&\disp\le\|v_1-v_2-2\kk_0(u_3-u_2)\|+2\kk_0\|u_3-u_1\|\\
&\disp=\|v_1-v_2-2\kk_0(\vt_0(v_1,p_2)-\vt_0(u_2,p_2)\|+2\kk_0\|\vt_0(v_1,p_2)-\vt_0(v_1,p_1)\|\\
&\disp\le\|v_1-v_2\|+2\kk_0\ell_0d(p_1,p_2)^\frac{1}{2},
\end{array}
\]
which show that the starting single-valued localization $\vt_0$ satisfies the H\"olderian inequality \eqref{Hol1}.

Now we are in a position to apply Lemma~\ref{lm1}(i) and to do propagation from $G_0$ to $G_\tau$ with $\tau=\frac{1}{2}$. Taking into account that the length of the propagation interval by Lemma~\ref{lm1} is $\th=\kk_0/(2\|B\|)$ and that the modulus $\ell$ in \eqref{Hol1} is doubled at each step, we need to make $n$ steps for reaching $G_{1/2}$ from $G_0$, where $n\in\IN$  is chosen from the interval
$$
\disp\frac{1}{2\th}\le n<\frac{1}{2\th}+1.
$$
In this way we get from Lemma~\ref{lm1}(i) that $G_\tau$ with $\tau=\frac{1}{2}$ admits a single-valued localization $\vt_\tau$ relative to a neighborhood $V_\tau\times Q_\tau\times U_\tau$ of $(\ov,\op,\ox)$ satisfying by \eqref{Hol2} the following inequality:
\begin{eqnarray}\label{4.35a}
\big\|v_1-v_2-2\kk_0\big[\vt_\tau(v_1,p_1)-\vt_\tau(v_2,p_2)\big]\big\|\le\|v_1-v_2\|+2^n(2\kk_0\ell_0)d(p_1p_2)^\frac{1}{2}
\end{eqnarray}
for all $(v_1,p_1),(v_2,p_2)\in V_\tau\times Q_\tau$. Finally, we employ Lemma~\ref{lm2}(i) to pass from the linearization $G_{1/2}$ to the solution map $S$ in \eqref{ss}. Since the modulus is $\kk_0$ on the left-hand side of \eqref{4.35a}, we choose any $\ve<\kk_0$ in \eqref{4.22a} and arrive at the H\"older property \eqref{4.7} with the assigned modulus $\kk=\kk_0-\ve$ in Definition~\ref{fs}(i), which can be chosen smaller than but arbitrarily close to $\kk_0$. This finally completes the proof of the theorem.\endproof

As a consequence of Theorem~\ref{Holder}, we get the equivalence between H\"olderian full stability of PVS \eqref{VS} when $f=0$ and the notion with the same name for local minimizers of $g$ defined above following \cite{MN}. In this way, the full stability conditions from Definition~\ref{fs}(i) can be seen as a new characterization of H\"olderian full stability of local minimizers in scalar optimization.

\begin{Corollary}{\bf (H\"olderian full stability of local minimizers).}\label{hol-min} The point $\ox$ is a H\"olderian fully stable local minimizer associated with $g\colon X\times P\to\oR$ relative to the parameter pair $(\ov,\op)$ if and only if it is a H\"olderian fully stable solution to PVS \eqref{VS} with $f=0$ corresponding to $(\ov,\op)$.
\end{Corollary}
{\bf Proof}. It follows from the comparison of the characterization of H\"olderian full stability in Theorem~\ref{Holder} with $f=0$ and that for H\"olderian full stability of local minimizers associated with $g$ obtained in \cite[Theorem~4.8]{MN} under the same assumptions on the function $g$ in \eqref{VS}.\endproof

Next we establish second-order characterizations of Lipschitzian full stability of solutions to PVS \eqref{VS} in the sense of Definition~\ref{fs}(ii). The first characterization concerns the general Hilbert space setting under the standing assumptions {\bf (A1)}--{\bf(A3)}. For brevity, the quantitative relationship between $\kk$ and $\kk_0$ similar to that of Theorem~\ref{Holder} is omitted in what follows.

\begin{Theorem}{\bf(neighborhood second-order characterization of Lipschitzian full stability of PVS).}\label{coro2} Let $\ox\in S(\ov,\op)$ be a solution to \eqref{VS} corresponding to the parameter pair $(\ov,\op)$. Then $\ox$ is Lipschitzian fully stable for \eqref{VS} if and only if the second-order condition \eqref{4.8} holds for some $\eta,\kk_0$  and in addition we have that the graphical partial subdifferential mapping
\begin{eqnarray}\label{K}
K\colon p\mapsto\gph\partial_x g(\cdot,p)\;\;\mbox{is Lipschitz-like around}\;\;(\op,\ox,\hat v),
\end{eqnarray}
where $\hat v=\ov-f(\ox,\op)$. Consequently, this implies that the solution map $S$ from \eqref{ss} admits a Lipschitz continuous and single-valued localization around $(\ov,\op,\ox)$.
\end{Theorem}
{\bf Proof.} To justify the ``only if" part of the theorem, it remains to show by the view of Theorem~\ref{Holder} that the Lipschitzian full stability of $\ox$ implies the Lipschitz-like property of the mapping $K$ in \eqref{K}. To proceed, fix the neighborhoods $U,Q,V$ from Definition~\ref{fs}(ii) and find a neighborhood  $U_1\times Q_1\times V_1\subset U\times Q\times X$ of $(\ox,\op,\hat v)$ such that $v+f(u,p)\in V$ for all $(u,p,v)\in U_1\times Q_1\times V_1$. Picking any $p_1,p_2\in Q_1$ and $(u_1,v_1)\in K(p_1)\cap(U_1\times V_1)$, it follows that
$$
V_1\ni v^\prime_1:=f(u_1,p_1)+v_1\in f(u_1,p_1)+\partial_x g(x_1,p_1)
$$
yielding $u_1=\vt(v_1^\prime,p_1)$. Define further $u_2:=\vt(v^\prime_1,p_2)$ and $v_2:=v^\prime_1-f(u_2,p_2)\in\partial_x g(u_2,p_2)$, which means $(u_2,v_2)\in K(p_2)$. Since $u_1=\vt(v_1^\prime,p_1)$ and $u_2=\vt(v^\prime_1,p_2)$, we deduce from \eqref{mp1} that
\begin{eqnarray}\label{XX}
2\kk\|u_1-u_2\|=\|(v^\prime_1-v^\prime_1)-2\kk(u_1-u_2)\|\le\|v^\prime_1-v^\prime_1\|+\ell d(p_1,p_2)=\ell d(p_1,p_2).
\end{eqnarray}
Moreover, it follows from the Lipschitzian property of $f$ in \eqref{4.6} that
\[
\|v_1-v_2\|=\|v_1-v^\prime_1+f(x_2,p_2)\|=\|-f(u_1,p_1)+f(u_2,p_2)\|\le L\big(\|u_1-u_2\|+d(p_1,p_2)\big)
\]
if $U_1,Q_1$ are chosen to be sufficiently small. This together with \eqref{XX} tells us that
\[
\|u_1-u_2\|+\|v_1-v_2\|\le (1+L)\|u_1-u_2\|+Ld(p_1,p_2)\le\frac{(1+L)\ell}{2\kk}d(p_1,p_2)+Ld(p_1,p_2).
\]
Thus we arrive at the inclusion
\begin{eqnarray*}
K(p_1)\cap (U_1\times V_1)\subset K(p_2)+\Big[\frac{(1+L)\ell}{2\kk}+L\Big]d(p_1,p_2)\B\quad\mbox{for all}\quad p_1,p_2\in Q_1,
\end{eqnarray*}
which verifies the claimed Lipschitz-like property of the mapping $K$ in \eqref{K}.\vspace*{0.03in}

Conversely, let us prove the sufficiency of conditions \eqref{4.8} and \eqref{K} for the Lipschitzian full stability of $\ox\in S(\ov,\op)$. It follows from Theorem~\ref{Holder} that there exist $\kk,\ell>0$ and a neighborhood $U\times Q\times V$ of $(\ox,\op,\ov)$ such that the H\"olderian condition \eqref{4.7}
is satisfied. Since $K$ is Lipschitz-like around $(\op,\ox,\hat v)$, we find $L_1>0$ and a neighborhood $U_1\times Q_1\times V_1\subset U\times Q\times X$ of $(\op,\ox,\hat v)$ for which $f(u,p)+v\in V$ whenever $(u,p,v)\in U_1\times Q_1\times V_1$ and
\begin{eqnarray}\label{4.39}
K(p_1)\cap(U_1\times V_1)\subset K(p_2)+L_1 d(p_1,p_2)\B\quad \mbox{for all}\quad p_1,p_2\in Q_1.
\end{eqnarray}
Take the localization $\vt$ from Definition~\ref{fs}(i) and get a neighborhood $U_2\times Q_2\times V_2\subset U_1\times Q_1\times V$ of $(\ox,\op,\ov)$ such that $\vt(V_2\times Q_2)\subset U_2$ and that $v-f(u,p)\subset V_1$ for all $(u,p,v)\in U_2\times Q_2\times V_2$. Pick now any $(v_1,p_1),(v_2,p_2)\in V_2\times Q_2$ and define $u_1:=\vt(v_1,p_1)\in U_2$ and $u_2:=\vt(v_2,p_2)\in U_2$. Therefore we have $v_1^\prime:=v_1-f(u_1,p_1)\in\partial_x g(u_1,p_1)$, i.e., $(u_1,v_1^\prime)\in K(p_1)$. It follows from \eqref{4.39} that there is $(u,v)\in K(p_2)$ satisfying the condition
\begin{eqnarray}\label{4.40}
\|u-u_1\|+\|v-v_1^\prime\|\le L_1d(p_1,p_2).
\end{eqnarray}
Define $v^\prime:=f(u,p_2)+v\in f(u,p_2)+\partial_x g(u,p_2)$ and observe from \eqref{4.6} and \eqref{4.40} that
\begin{eqnarray}\label{VU}\begin{array}{ll}
\|v^\prime-v_1\|&\disp=\|f(u,p_2)+v-f(u_1,p_1)-v_1^\prime\|\le\|v-v^\prime_1\|+\|f(u,p_2)-f(u_1,p_1)\|\\
&\disp\le L_1d(p_1,p_2)+L\big(\|u-u_1\|+d(p_1,p_2)\big)\\
&\disp\le L_1d(p_1,p_2)+L\big(L_1d(p_1,p_2)+d(p_1,p_2)\big)\\
&\disp=(L_1+LL_1+L)d(p_1,p_2).
\end{array}
\end{eqnarray}
Hence we have $v^\prime\in V$ by choosing $Q_2$ to be sufficiently small, which implies that $u=\vt(v^\prime,p_2)$. Using now \eqref{4.7} for the pairs $(v',p_2)$ and $(v_2,p_2)$ gives us the inequality
\begin{eqnarray*}
\|(v^\prime-v_2)-2\kk(u-u_2)\|\le\|v^\prime-v_2\|.
\end{eqnarray*}
Combining this with \eqref{VU} and \eqref{4.40} ensures the estimates
\begin{eqnarray*}\begin{array}{ll}
\|(v_1-v_2)-2\kk(u_1-u_2)\|&\disp\le\|(v^\prime-v_2)-2\kk(u-u_2)\|+\|v^\prime-v_1\|+2\kk\|u_1-u\|\\
&\disp\le\|v^\prime -v_2\|+\|v^\prime-v_1\|+2\kk\|u_1-u\|\\
&\disp\le\|v_1-v_2\|+2\|v^\prime-v_1\|+2\kk\|u_1-u\|\\
&\disp\le\|v_1-v_2\|+2(L_1+LL_1+L)d(p_1,p_2)+2\kk L_1 d(p_1,p_2),
\end{array}
\end{eqnarray*}
which verify \eqref{mp1} and thus complete the proof of the theorem.\endproof

Similarly to Corollary~\ref{hol-min}, we can establish the equivalence between the Lipschitzian full stability of solutions to PVS \eqref{VS} with $f=0$ and Lipschitzian full stability of local minimizers \eqref{LS} for the corresponding optimization problem associated with $g$. This follows from the comparison of the second-order characterizations obtained in Theorem~\ref{coro2} and in \cite[Corollary~4.8]{MN}, respectively.\vspace*{0.05in}

Observe by the coderivative criterion \eqref{cod-cr} that in the case of both finite-dimensional spaces $X$ and $P$ the Lipschitz-like property \eqref{K} is equivalent to the pointwise second-order condition
\begin{eqnarray}\label{Mor}
(0,q)\in\big(D^*\partial_x g\big)(\ox,\op,\hat v)(0)\Longrightarrow q=0.
\end{eqnarray}

Our next result provides a complete pointwise characterization of Lipschitzian full stability of PVS via the limiting second-order subdifferential constructions for $g$.

\begin{Theorem}{\bf (pointwise characterization of Lipschitzian full stability of PVS in finite dimensions).} \label{Lips} Let $X,P$ be two finite-dimensional spaces, and let $\ox\in S(\ov,\op)$. Then the Lipschitzian full stability of $\ox$ is equivalent to the simultaneous validity of \eqref{Mor} and the condition
\begin{eqnarray}\label{4.43}
\la\nabla_x f(\ox,\op)w,w\ra+\la z,w\ra>0\quad\mbox{for all}\quad(z,q)\in\big(D^*\partial_x g\big)(\ox,\op,\hat v)(w),\;w\ne 0.
\end{eqnarray}
Consequently, conditions \eqref{Mor} and \eqref{4.43} imply that the solution map $S$ from \eqref{ss} admits a Lipschitz continuous and single-valued localization around $(\ov,\op,\ox)$.
\end{Theorem}
{\bf Proof.} Similarly to the proof of Theorem~\ref{Holder}, consider the l.s.c.\ function $h\colon X\times P\to\oR$ defined in \eqref{hh} and easily get from the elementary sum rule for the subdifferential \eqref{2.2} that
$$
\partial_x h(x,p)=f(\ox,\op)+\frac{1}{2}\Big(A+A^*\Big)(x-\ox)+\partial_x g(x,p)\;\mbox{ and }\;\ov\in\partial_x h(\ox,\op).
$$
As mentioned above, $h$ is parametrically continuously prox-regular at $(\ox,\op)$ for $\ov$, and it follows from \cite[Theorem~1.62]{M1} the representation
\begin{eqnarray}\label{cod}
\big(D^*\partial_x h\big)(\ox,\op,\ov)(w)=\Big(\frac{1}{2}(A+A^*)w,0\Big)+\big(D^*\partial_x g\big)(\ox,\op,\hat v)(w).
\end{eqnarray}

To justify the ``only if" part of the theorem, suppose that $\ox$ is a Lipschitzian full stable solution to \eqref{VS} corresponding to $(\ov,\op)$ and find a single-valued localization $\vt$ of $S$ satisfying \eqref{mp1}. It follows from Theorem~\ref{coro2} that condition \eqref{K} holds and so does \eqref{Mor}. By \eqref{cod} we easily have
\begin{eqnarray}\label{Mor2}
(0,q)\in\big(D^*\partial_x h\big)(\ox,\op,\ov)(0)\Longrightarrow q=0.
\end{eqnarray}
Note also that the Lipschitzian full stability of $\ox$ ensures the validity of condition \eqref{4.8} in Theorem~\ref{Holder}. Furthermore, it is shown in the proof of Theorem~\ref{Holder} that the second-order growth condition \eqref{4.33} holds. When the latter two conditions are satisfied, we have that $\ox$ is a Lipschitzian fully stable local minimizer of $h$ relative to $\op$ and $\ov$ by \cite[Theorem~3.4 and Theorem~4.6]{MN}. Then employing \cite[Corollary~4.10]{MN} for the function $h$ gives us
\begin{eqnarray}\label{4.46a}
\la z,w\ra>0\quad\mbox{for all}\quad(z,q)\in\big(D^*\partial_x h\big)(\ox,\op,\ov)(w),\;w\ne 0.
\end{eqnarray}
This together with equality \eqref{cod} implies that
\begin{eqnarray*}
0<\frac{1}{2}\la(A+A^*)w,w\ra+\la z,w\ra=\la A w,w\ra+\la z,w\ra\;\mbox{for all}\;z\in\big(D^*\partial_x g\big)(\ox,\op,\hat v)(w),\;w\ne 0,
\end{eqnarray*}
which ensures \eqref{4.43} and thus completes the proof of the necessity part of the theorem.

To verify the sufficiency of conditions \eqref{Mor} and \eqref{4.43} for the Lipschitzian full stability of $\ox$, suppose that these conditions are satisfied and mention again that \eqref{Mor} implies \eqref{Mor2} due to equality \eqref{cod}. Similarly to the above but in the opposite direction, we can verify the validity of \eqref{4.46a} from \eqref{4.43} and \eqref{cod}. Now we can use the characterizations of Lipschitzian first stability of local minimizers associated with $h$, first from \cite[Corollary~4.10 and Theorem~4.6]{MN} and then from \cite[Theorem~3.4]{MN}, which allow us to conclude that the mapping $G_0$ in \eqref{gt} admits a single-valued and Lipschitz continuous localization $\vt_0$ relative to a neighborhood $V\times Q\times U$ of $(\ov,\op,\ox)$ with some modulus $\kk_0>0$ such that for any triple $(u,p,v)\in\gph\vt_0$ we have
\begin{eqnarray*}
h(x,p)\ge h(u,p)+\la v,x-u\ra+\frac{\kk_0}{2}\|x-u\|^2\quad\mbox{whenever}\quad x\in U.
\end{eqnarray*}
The latter implies similarly to the proof of \eqref{m2} the strong monotonicity condition
\begin{eqnarray*}
\la v_1-v_2,u_1-u_2\ra\ge\kk_0\|u_1-u_2\|^2,
\end{eqnarray*}
which ensures as in the subsequent proof of Theorem~\ref{Holder} that the mapping $G_0$ admits a single-valued localization $\vt_0$ satisfying condition \eqref{Lip1}. Employing in this vein the propagation from Lemma~\ref{lm1}(ii) finitely many times shows that the mapping $G_\tau$ with $\tau=\frac{1}{2}$ also has a single-valued localization $\vt_\tau$ satisfying \eqref{Lip1}. Finally, Lemma~\ref{lm2}(ii) allows us to pass to the solution map $S$ in \eqref{ss} and completes the proof of the theorem.\endproof

When the parameter $p$ is omitted in \eqref{VS}, we derive from Theorem~\ref{Lips} the following new characterization of local strong maximal monotonicity for the important class of set-valued mappings.

\begin{Corollary}{\bf(characterization of local strong maximal monotonicity of set-valued mappings).}\label{coro4} Let $X$ be a finite-dimensional space, let $f:X\to X$ be continuously differentiable around $\ox$, and let $g\colon X\to\oR$ be prox-regular and subdifferentially continuous at $\ox$ for some vector $\hat v\in\partial g(\ox)$. The following assertions are equivalent:

{\bf (i)} The set-valued map $T=f+\partial g:X\tto X$ is locally strongly maximally monotone around $(\ox,\ov)$ with $\ov:=f(\ox)+\hat v$.

{\bf (ii)} The limiting coderivative $D^*T(\ox,\ov)$ is positive-definite in the sense of \eqref{pd}.
\end{Corollary}
{\bf Proof.} It follows from Lemma~\ref{pro2} that assertion {\bf (i)} of this corollary is equivalent to condition \eqref{3.5}, which exactly is the Lipschitzian full stability condition from Definition~\ref{fs}(ii) in the nonparametric case. Theorem~\ref{Lips} tells us that the latter is equivalent to
\[
\la\nabla f(\ox)w,w\ra+\la z,w\ra>0\quad\mbox{for all}\quad z\in\big(D^*\partial g\big)(\ox,\hat v)(w),\;w\ne 0,
\]
which reduces to \eqref{pd} for $T=f+\partial g$ by the coderivative sum rule in \cite[Theorem~1.62]{M1}.\endproof

Note that Corollary~\ref{coro4} provides yet another evidence of the validity of Conjecture~\ref{con} formulated above. To see it, we need only observing that that the mapping $T=f+\partial g$ is hypomonotone around $(\ox,\ov)$ due to the hypomonotonicity of limiting subgradient mappings for continuously prox-regular functions and the hypomonotonicity sum rule from Proposition~\ref{pro1}(ii).

\section{Full Stability in Parametric Variational Inequalities}
\setcounter{equation}{0}

In this section we consider a particular case of the general parametric variational systems \eqref{VS} described by {\em parametric variational inequalities} (PVI) in the form
\begin{eqnarray}\label{Ge}
v\in f(x,p)+N_C(x),
\end{eqnarray}
where $X$ is a Hilbert space and $P$ is a metric space endowed with some metric $d$. Besides this, our {\em standing assumptions} here are that
the mapping $f:X\times P\to X$ satisfies {\bf (A1)} from Section~4 and that $C$ is a closed and {\em convex} subset of $X$. As discussed in Section~1, the variational system (generalized equation) \eqref{Ge} can be rewritten in the standard variational inequality form \eqref{vi1}.

Denote the solution map to the parametric variational inequality \eqref{Ge} by
\begin{eqnarray}\label{sol}
\Hat S(v,p):=\big\{x\in X\big|\;v\in f(x,p)+N_C(x)\big\}.
\end{eqnarray}
The main goal of this section is to show that necessary and sufficient conditions for Lipschitzian full stability of solutions $\ox\in\hat S(\ov,\op)$ to \eqref{Ge} can be obtained in the {\em pointwise} form as in \eqref{4.43} in {\em infinite-dimensional} decision spaces $X$ under some reasonable assumptions, which have been well understood and applied in the case of infinite-dimensional variational and control problems, especially those related to semilinear partial differential equations of the elliptic type. Furthermore, we establish conditions for Lipschitzian full stability in the PVI setting under consideration in more explicit forms involving the underlying convex set $C$ generated the variational inequality \eqref{Ge}.

Since in this and next sections we focus only on {\em Lipschitzian} full stability of the corresponding parametric variational systems, we will omit the word ``Lipschitzian" in what follows.\vspace*{0.01in}

To proceed, we first recall the aforementioned assumptions following \cite[p.\ 194]{BS} and \cite[p.\ 259]{IT} for the Legendre form and \cite{H,Mi} for polyhedric sets; see also \cite{BBS,B,BS,hms,HS,MN} and the references therein for more details, discussions, and applications.

\begin{Definition}{\bf (Legendre forms).}\label{legen} The real-valued function $Q:X\to\R$ is a {\sc Legendre form} if it is weakly lower semicontinuous, represented as $Q(x)=\la Ax,x\ra$ with some linear operator $A:X\to X$, and satisfies the implication
$$
\big[x_k\st{w}\to x,\;Q(x_k)\to Q(x)\big]\Longrightarrow x_k\to x\;\mbox{ as }\;k\to\infty.
$$
\end{Definition}

We see that the conditions of Definition~\ref{legen} hold trivially for quadratic forms in finite dimensions. In general Hilbert spaces, various necessary and sufficient conditions for the validity of Definition~\ref{legen} and interesting examples of Legendre forms can be found in \cite[Section~3.3.2]{BS} and \cite[Section~6.2]{IT}.

\begin{Definition} {\bf (polyhedric sets).}\label{poly} Let $C$ be a closed and convex subset of $X$. It is said to be {\sc polyhedric} at $\ox\in C$ for some $\hat v\in N_C(\ox)$ if we have the representation
\begin{eqnarray}\label{5.3}
\K_C(\ox,\hat v):=T_C(\ox)\cap\{\hat v\}^\perp={\rm cl}\Big\{\mathcal{R}_C(\ox)\cap\{\hat v\}^\perp\Big\}
\end{eqnarray}
of the corresponding critical cone $\K(\ox,\hat v)$, where
\begin{equation}\label{5.3a}
\mathcal{R}_C(\ox):=\bigcup_{t>0}\Big[\frac{C-\ox}{t}\Big]
\end{equation}
is called the radial cone, and where $T_C(\ox):={\rm cl}\,\mathcal{R}_C(\ox)$ is the classical tangent cone to $C$ at $\ox$. If $C$ is polyhedric at each $\ox\in C$ for any $\hat v\in N_C(\ox)$, we say that the set $C$ is polyhedric.
\end{Definition}

Polyhedral sets in finite and infinite dimensions are automatically polyhedric while the latter class is significantly broader; see, in particular, \cite[Chapter~6]{BS} and \cite{BBS,B,H,HS,Mi} for important examples of polyhedric but nonpolyhedral sets in infinite dimensions.\vspace*{0.03in}

The following proposition, taken from \cite[Theorem~6.2]{MN} and employed below, provides a precise calculation of the regular coderivative for the normal cone mapping generated by a polyhedric set.

\begin{Proposition}{\bf (regular coderivative calculation for the normal cone mapping to polyhedric sets).}\label{pro2.3} For any $\ox\in C$ and $\hat v\in N_C(\ox)$ we have the inclusion
\begin{eqnarray}\label{5.4}
\dom\hat D^*N_C(\ox,\hat v)\subset-\K_C(\ox,\hat v).
\end{eqnarray}
If in addition $C$ is polyhedric at $\ox\in C$ for $\hat v$, then
\begin{eqnarray}\label{5.5}
\Hat D^*N_C(\ox,\hat v)(w)=\K_C(\ox,\hat v)^*\quad \mbox{whenever}\quad w\in -\K_C(\ox,\hat v).
\end{eqnarray}
\end{Proposition}

Now we are in a position to derive a major result of this section that gives pointwise sufficient and necessary conditions for full stability of solutions to PVI \eqref{Ge} in Hilbert spaces. It can be seen as a far-going extension of \cite[Theorem~7.2]{MN}, which characterizes (Lipschitzian) full stability of local solutions in the sense of \eqref{LS}, to optimal control problems governed by semilinear partial differential equations with elliptic operators written as a Legendre form.

\begin{Theorem}{\bf (pointwise sufficient and necessary conditions for full stability of solutions to PVI).}\label{NGE} Let $\ox\in\Hat S(\ov,\op)$ in \eqref{sol}, and let $\Hat v:=\ox-f(\ox,\op$. Consider the following statements:

{\bf (i)} $\ox$ is a fully stable solution to PVI \eqref{Ge} in the sense of Definition~{\rm\ref{fs}}{\rm(ii)}.

{\bf (ii)} We have the positive-definiteness condition
\begin{eqnarray}\label{5.7}
\la\nabla_x f(\ox,\op)w,w\ra>0\quad\mbox{for all}\quad w\in{\cal H}^{\rm w}(\ox,\hat v),\;w\ne 0,
\end{eqnarray}
in terms of the $($weak$)$ outer limit \eqref{pk} of the critical cones
\begin{eqnarray}\label{5.8}
{\cal H}^{\rm w}(\ox,\hat v):=\Limsup_{(x,v)\st{{\rm gph}\,N_C}\longrightarrow(\ox,\hat v)}\K_C(x,v).
\end{eqnarray}

{\bf (iii)} We have the  positive-definiteness condition
\begin{eqnarray}\label{5.9}
\la\nabla_x f(\ox,\op)w,w\ra>0\quad \mbox{for all}\quad w\in{\cal H}^{\rm s}(\ox,\hat v),\;w\ne 0,
\end{eqnarray}
where ${\cal H}^{\rm s}(\ox,\hat v)$ is defined by the the strong version of the outer limit \eqref{pk} via the norm topology
\begin{eqnarray}\label{5.10}
\begin{array}{ll}
{\cal H}^{\rm s}(\ox,\hat v)&\disp:={\rm s}-\Limsup_{(x,v)\st{{\rm gph}\,N_C}\longrightarrow(\ox,\hat v)}\K_C(x,v)\\
&\disp=\Big\{z\in X\big|\;\exists\;{\rm sequences }\,(x_k,v_k)\st{{\rm gph}N_C}\longrightarrow(\ox,\hat v),\; z_k\in \K_C(x_k,v_k),\;z_k\to z\Big\}.
\end{array}
\end{eqnarray}
Then the following assertions are satisfied:

{\bf(1)} If $Q(w):=\la\nabla_x f(\ox,\op)w,w\ra$ is Legendre, then {\bf (ii)} is a sufficient condition for {\bf (i)}.

{\bf(2)} If  $C$ is polyhedric, then  {\bf (iii)} is necessary for the validity of {\bf (i)}.
\end{Theorem}
{\bf Proof.} To justify assertion~{\bf(1)}, let $Q$ be Legendre, and let \eqref{5.7} hold. Since BCQ \eqref{bcq} is trivial when $g(x,p):=\delta_C(x)$ and so $\partial_x g(x,p)=N_C(x)$ for all $(x,p)\in X\times P$, it follows from Theorem~\ref{coro2} that we only need to verify the validity of \eqref{4.8}. Assume the contrary and find a sequence $(u_k,v_k,w_k,z_k)\in X\times X\times X\times X$ such that $(u_k,v_k)\st{{\rm gph}\,N_C}\to (\ox,\hat v)$, $z_k\in\Hat D^*N_C(u_k,v_k)(w_k)$, and
\begin{eqnarray}\label{5.11}
\la\nabla_x f(\ox,\op)w_k,w_k\ra+\la z_k,w_k\ra<\frac{1}{k}\|w_k\|^2\;\mbox{ for all }\;k\in\IN.
\end{eqnarray}
Since $N_C$ is maximally monotone, we get from \cite[Lemma~3.3]{CT} that $\la z_k,w_k\ra\ge 0$. Furthermore, it follows from Proposition~\ref{pro2.3} that $w_k\in-{\cal K}_C(u_k,v_k)$. Combining these facts with \eqref{5.11} yields
\begin{eqnarray}\label{Q0}
\la\nabla_x f(\ox,\op)w_k,w_k\ra<\frac{1}{k}\|w_k\|^2\quad\mbox{with}\quad w_k\in-{\cal K}_C(u_k,v_k),
\end{eqnarray}
which readily implies that $w_k\ne 0$. Defining $\ow_k:=w_k\|w_k\|^{-1}$, we get from \eqref{Q0} that
\begin{eqnarray}\label{Q}
Q(\ow_k)=\la\nabla_x f(\ox,\op)\ow_k,\ow_k\ra<\frac{1}{k}\quad\mbox{with}\quad\ow_k\in-{\cal K}_C(u_k,v_k).
\end{eqnarray}
Since $\|\ow_k\|=1$, there is a subsequence of $\{\ow_k\}$ (no relabeling), which weakly converges to some $\ow$. By $(u_k,v_k)\st{{\rm gph}\, N_C}\to(\ox,\hat v)$, we derive from \eqref{5.8} that $\ow\in-{\cal H}^{\rm w}(\ox,\hat v)$. Employing now the weak l.s.c.\ of the Legendre form $Q$ ensures by \eqref{Q} and \eqref{5.7} that
\begin{eqnarray}\label{Q1}
0\le Q(\ow)\le\liminf_{k\to\infty}Q(\ow_k)\le 0,
\end{eqnarray}
which yields $Q(\ow)=0$. By passing to a subsequence in \eqref{Q1} and Definition~\ref{legen} for $Q$, suppose without loss of generality that $\ow_k\to\ow$ as $k\to\infty$. Hence we get $\|\ow\|=1$, $\ow\in{\cal H}^{\rm w}(\ox,\hat v)$, and $Q(\ow)=Q(-\ow)=0$. This contradicts \eqref{5.7} and thus completes the proof of assertion~{\bf(1)}.\vspace*{0.03in}

It remains to verify assertion~{\bf(2)} provided that $C$ is polyhedric. Indeed, by the assumed full stability in {\bf (i)} it follows from Theorem~\ref{coro2} in the case of $g(x,p)=\delta_C(x)$ that there are numbers $\kk,\eta>0$ such that for any $(u,v)\in\gph N_C\cap\B_\eta(\ox,\hat v)$ we have
\begin{eqnarray}\label{5.13}
\la\nabla_x f(\ox,\op)w,w\ra+\la z,w\ra\ge\kk\|w\|^2\quad\mbox{whenever}\quad z\in\Hat D^*N_C(u,v)(w),w\in X.
\end{eqnarray}
Since $C$ is polyhedric, Proposition~\ref{pro2.3} tells us $w\in-{\cal K}(u,v)$ and $\Hat D^*N_C(u,v)(w)={\cal K}_C(u,v)^*$, which yield $0\in\Hat D^*N_C(u,v)(w)$. This together with \eqref{5.13} and \eqref{5.5} shows that
\begin{eqnarray}\label{5.13a}
\la\nabla_x f(\ox,\op)w,w\ra\ge\kk\|w\|^2\quad\mbox{for all}\quad w\in-{\cal K}_C(u,v)
\end{eqnarray}
for all $(u,v)\in\gph N_C\cap\B_\eta(\ox,\hat v)$. Employing the limiting procedure in \eqref{5.13a} as $\eta\dn 0$ and using the strong convergence in the definition of ${\cal H}^{\rm s}$, we arrive at \eqref{5.9} and so complete the proof.\endproof

The next lemma effectively estimates the limiting forms ${\cal H}^{\rm w}$ and ${\cal H}^{\rm s}$ in Theorem~\ref{NGE} via the tangent and critical cones for the set $C$, which allows us to establish more direct and verifiable conditions for full stability of solutions to \eqref{Ge} in both finite and infinite dimensions.

\begin{Lemma}{\bf(estimates of weak and strong outer limits of the critical cone).}\label{pro5.5} In the general setting of Theorem~{\rm\ref{NGE}} we have the upper estimate
\begin{eqnarray}\label{sub}
{\cal H}^{\rm w}(\ox,\hat v)\subset{\rm cl}\big[T_C(\ox)-T_C(\ox)\big]\cap\{\hat v\}^\perp.
\end{eqnarray}
If in addition $C$ is polyhedric, then the lower estimate
\begin{eqnarray}\label{sup}
{\cal K}_C(\ox,\hat v)-{\cal K}_C(\ox,\hat v)\subset{\cal H}^{\rm s}(\ox,\hat v)
\end{eqnarray}
is satisfied. It furthermore $\dim X<\infty$ and the set $C$ is polyhedral, then we have the representation
\begin{eqnarray}\label{5.16}
{\cal H}^{\rm s}(\ox,\hat v)={\cal H}^{\rm w}(\ox,\hat v)={\cal K}_C(\ox,\hat v)-{\cal K}_C(\ox,\hat v).
\end{eqnarray}
\end{Lemma}
{\bf Proof.} To verify \eqref{sub}, pick any $w\in{\cal H}^{\rm w}(\ox,\hat v)$ and find a sequence $(u_k,v_k,w_k)$ such that $w_k\st{w}\to w$, $(u_k,v_k)\st{{\rm gph}N_C}\to(\ox,\hat v)$, and $w_k\in{\cal K}_C(u_k,v_k)=T_C(u_k)\cap\{v_k\}^\perp$. Hence for each $k\in\IN$ there exist sequences $y_{n_k}\to u_k$ and $t_{n_k}\dn 0$ with $\frac{y_{n_k}-u_k}{t_{n_k}}\to w_k$ as $n\to\infty$. In this way we construct sequences
$$
(z_k,\al_k)\in\big\{(y_{n_k},t_{n_k})\big|\;n\in\IN\big\}
$$
so that $z_k\to\ox$, $\al_k\dn 0$ , and $\frac{z_{k}-u_k}{\al_k}-w_k\to 0$ as $k\to\infty$. It follows that
\[
w={\rm w}-\lim_{k\to\infty}w_k={\rm w}-\lim_{k\to\infty}\frac{z_{k}-u_k}{\al_k}={\rm w}-\lim_{k\to\infty}\frac{(z_{k}-\ox)-(u_k-\ox)}{\al_k}\in {\rm cl}^{\rm w}\,\Big[T_C(\ox)-T_C(\ox)\Big],
\]
where the symbol ``${\rm w}-\lim$" indicates that the weak limit in $X$ is taken. Since $\la w_k,v_k\ra=0$, we get $\la w,\hat v\ra=0$ by passing to the limit, and hence
$$
w\in {\rm cl}^{\rm w}\big[T_C(\ox)-T_C(\ox)\big]\cap\{\hat v\}^\perp.
$$
Note that $T_C(\ox)-T_C(\ox)$ is a convex set. The classical Mazur theorem tells us that ${\rm cl}^{\rm w}\big[T_C(\ox)-T_C(\ox)\big]={\rm cl}\big[T_C(\ox)-T_C(\ox)\big]$. Therefore $w\in{\rm cl}\big[T_C(\ox)-T_C(\ox)\big]\cap\{\hat v\}^\perp$, which justifies \eqref{sub}.\vspace*{0.03in}

Now suppose that $C$ is polyhedric. To verify \eqref{sup}, pick any $w=w_1-w_2$ from the left-hand side of \eqref{sup} with $w_1,w_2\in {\cal K}_C(\ox,\hat v)$. The polyhedricity of $C$ allows us to find \eqref{5.3} sequences $w_{1k}\to w_1$, $w_{2k}\to w_2$, and $t_{1k},t_{2k}\dn 0$ such that $\ox+t_{1k}w_{1k}\in C$, $\ox+t_{2k}w_{2k}\in C$, and $w_{1k}, w_{2k}\in\{\hat v\}^\perp$. Defining $t_k:=\min\{t_{1k},t_{2k}\}$, we get from the convexity of $C$ that $u_k:=\ox+t_kw_{1k}=(1-t_kt_{1k}^{-1})\ox+t_kt_{1k}^{-1}(\ox+t_{1k}w_{1k})\in C$ and similarly $x_k:=\ox+t_kw_{2k}\in C$. Thus it follows from \eqref{5.3} and \eqref{5.3a} that
\[
w_{2k}-w_{1k}=\frac{x_k-u_k}{t_k}\in{\cal R}_C(u_k)\cap\{\hat v\}^\perp\subset\K_C(u_k,\hat v).
\]
Since $w_{2k}-w_{1k}\to w_2-w_1=w$ and $u_k\to\ox$ as $k\to\infty$, we deduce from \eqref{5.10} that $w\in{\cal H}^{\rm s}(\ox,\hat v)$ and hence complete the proof of inclusion \eqref{sup}.

Finally, assuming $\dim X<\infty$ ensures that ${\cal H}^{\rm w}(\ox,\hat v)={\cal H}^{\rm s}(\ox,\hat v)$. To verify \eqref{5.16}, it remains to show  by \eqref{sup} that the opposite inclusion holds therein when $C$ is polyhedral. Picking any $w\in{\cal H}^{\rm s}(\ox,\hat v)$, find a sequence $(w_k, x_k,v_k)$ such that $w_k\in{\cal K}_C(x_k,v_k)$, $(x_k,v_k)\st{{\rm gph}N_C}\to(\ox,\hat v)$ and $w_k\to w$ as $k\to\infty$. It follows from the proof of \cite[Theorem~2]{DR1} that
\begin{eqnarray}
&&w_k\in T_C(x_k)=T_C(\ox)+\R\{x_k-\ox\},\label{TC}\\
&&v_k\in N_C(x_k)=N_C(\ox)\cap\{x_k-\ox\}^\perp\label{NC}
\end{eqnarray}
for all large $k\in\IN$ and that $T_C(\ox)\cap\{v_k\}^\perp\subset T_C(\ox)\cap\{\hat v\}^\perp={\cal K}_C(\ox,\hat v)$. Since $\R_+\{x_k-\ox\}\subset T_C(\ox)$ and the cone $T_C(\ox)$ is convex, we obtain from \eqref{TC} that
\[
w_k\in T_C(x_k)=T_C(\ox)-\R_+\{x_k-\ox\}.
\]
Hence there exist $y_k\in T_C(\ox)$ and $r_k\ge 0$ with $w_k=y_k-r_k(x_k-\ox)$. Since $v_k\in N_C(x_k)$, we get $\la x_k-\ox,v_k\ra=0$ from \eqref{NC}. It follows furthermore that $0=\la w_k,v_k\ra=\la y_k,v_k\ra$, and so $y_k\in T_C(\ox)\cap\{v_k\}^\perp\subset{\cal K}_C(\ox,\hat v)$. Observing from \eqref{NC} that $r_k(x_k-\ox)\subset T_C(\ox)\cap\{v_k\}^\perp\subset{\cal K}(\ox,\hat v)$ and using the above representation of $w_k$ give us the inclusion
\[
w_k\in{\cal K}_C(\ox,\hat v)-{\cal K}_C(\ox,\hat v)
\]
showing that $w\in{\rm cl}\,\big[{\cal K}_C(\ox,\hat v)-{\cal K}_C(\ox,\hat v)\big]={\cal K}_C(\ox,\hat v)-{\cal K}_C(\ox,\hat v)$, where the equality holds since ${\cal K}_C(\ox,\hat v)-{\cal K}_C(\ox,\hat v)$ is a subspace in finite dimensions. This verifies that\ \eqref{5.16} is satisfied.\endproof

The obtained results lead us to the following verifiable conditions for full stability in \eqref{Ge}.

\begin{Theorem}{\bf(refined pointwise conditions for full stability of solutions to PVS).}\label{coro5.6} In the framework of Theorem~{\rm\ref{NGE}}, consider the statements:

{\bf (i)} $\ox\in\Hat S(\ov,\op)$ is a fully stable solution to \eqref{Ge}.

{\bf (ii)} We have the positive-definiteness condition of closure type
\begin{eqnarray}\label{5.19}
\la\nabla_x f(\ox,\op)w,w\ra>0\quad\mbox{for all}\quad w\in{\rm cl}\big[T_C(\ox)-T_C(\ox)\big]\cap\{\hat v\}^\perp,\; w\ne 0.
\end{eqnarray}

{\bf (iii)} We have the the positive-definiteness condition via the critical cone
\begin{eqnarray}\label{5.20}
\la\nabla_x f(\ox,\op)w,w\ra>0\quad\mbox{for all}\quad w\in{\cal K}_C(\ox,\hat v)-{\cal K}_C(\ox,\hat v),\;w\ne 0.
\end{eqnarray}
Then the following assertions hold:

{\bf(1)} If $Q(w)=\la\nabla_x f(\ox,\op)w,w\ra$ is Legendre, then {\bf (ii)} is sufficient for the validity of {\bf(i)}.

{\bf(2)} If $C$ is polyhedric, then {\bf (iii)} is necessary for the validity of {\bf (i)}.

{\bf(3)} If $\dim X<\infty$ and $C$ is polyhedral, then the conditions in {\bf(ii)} and {\bf (iii)} are equivalent being necessary and sufficient for the validity of {\bf (i)}.
\end{Theorem}
{\bf Proof.} It follows by combining the corresponding assertions in Theorem~\ref{NGE} and Lemma~\ref{pro5.5}.\endproof

Note that the positive-definiteness condition of the closure type \eqref{5.19} has been used in \cite[p.\ 405]{BS} to ensure that the solution map $\Hat S$ in \eqref{sol} admits a single-valued and Lipschitz continuous localization around $(\ov,\op)$ provided that $Q(w)$ is Legendre. As follows from the discussion in Section~4, our assertion {\bf (1)} of Theorem~\ref{coro5.6} establishes much stronger property of {\em full stability} of solutions to \eqref{Ge} under the same assumptions as in \cite{BS}.\vspace*{0.03in}

Finally in this section, we focus on the space $X=L^2(\O)$, where $\O$ is an open subset of $\R^n$. Consider the closed and convex set $C\subset X$ in \eqref{Ge} defined by the {\em magnitude constraint system}
\begin{eqnarray}\label{C}
C=\big\{x\in L^2(\O)\big|\;a\le x(y)\le b\;\mbox{ a.e. }\;y\in\O\big\}
\end{eqnarray}
with $-\infty\le a<b\le\infty$. Such constraints are typical in applications to PDE control of elliptic equations; see, e.g., \cite{B,BS,MN}.
It follows from \cite[Proposition~6.33]{BS} that the set $C\subset L^2(\O)$ in \eqref{C} is polyhedric. The next proposition taken from \cite[Proposition~7.3]{MN} gives us precise calculations of the strong and weak outer limits of critical cones in \eqref{5.8} and \eqref{5.10}, respectively.

\begin{Proposition}{\bf(calculations of outer limits of critical cones).}\label{lim-calc} Let $(\ox,\hat v)\in\gph N_C$ with $C$ defined in \eqref{C}. Then we have the set \eqref{5.8} and \eqref{5.10} are calculated by
\begin{eqnarray*}
{\cal H}^{\rm s}(\ox,\hat v)={\cal H}^{\rm w}(\ox,\hat v)=\big\{u\in L^2(\O)\big|\; u(y)\hat v(y)=0\;\;\mbox{a.e.}\;\; \mbox{on}\;\;\O\big\}.
\end{eqnarray*}
\end{Proposition}

This allows us to derive from Theorem~\ref{NGE} the following pointwise characterization of full stability of solutions to infinite-dimensional PVS \eqref{Ge} generated by the constrained sets of type \eqref{C}.

\begin{Corollary}{\bf (characterization of full stability for PVI generated by magnitude constrained systems).} Let $X=L^2(\O)$ with $C$ given in  \eqref{C}, and let $\ox\in\hat S(\ov,\op)$ in \eqref{sol}. In addition to the standing assumption {\bf(A1)} on $f\colon X\times P\to X$ with a metric parameter space $P$, suppose that $Q(w):=\la\nabla f(\ox,\op)w,w\ra$ for $w\in L^2(\O)$ is a Legendre form. Then $\ox$ is a fully stable solution to \eqref{Ge} if and only if we have the pointwise positive-definiteness condition
\begin{eqnarray*}
\la\nabla_x f(\ox,\op)w,w\ra>0\quad\mbox{for all}\quad w\in L^2(\O)\setminus\{0\}\quad\mbox{with}\quad w(y)\big(\ov(y)-f(\ox,\op)(y)\big)=0\;\; \mbox{a.e.}\;\;y\in\O.
\end{eqnarray*}
\end{Corollary}
{\bf Proof.} Follows directly from Theorem~\ref{NGE} by employing the calculations of Proposition~\ref{lim-calc}.\endproof

\section{Full Stability in Parametric Variational Conditions}
\setcounter{equation}{0}

Here we study the {\em parametric variational conditions} (PVC) given in the form
\begin{eqnarray}\label{VC}
v\in f(x,p)+N_{C(p)}(x),
\end{eqnarray}
where $x\in X=\R^n$, $p\in P=\R^d$, and the set $C(p)$ for each $p$ is defined by the inequality constraints
\begin{eqnarray}\label{6.2}
C(p):=\big\{x\in\R^n\big|\;\ph_i(x,p)\le 0\;\mbox{ for }\;i=1,\ldots,m\big\}
\end{eqnarray}
via the functions $\ph_i:\R^n\times\R^d\to\R$, $i=1,\ldots,m$, which are ${\cal C}^2$-smooth around a feasible point $(\ox,\op)\in\R^n\times\R^d$ of \eqref{6.2}. As mentioned above, inclusion \eqref{VC} can be written in form \eqref{VS} with $g(x,p)=\delta_{C(p)}(x)$ for $(x,p)\in\R^n\times\R^d$. Furthermore, \eqref{VC} can be represented as a quasi-variational inequality when the sets $C(p)$ are convex. The solution map to \eqref{VC} is denoted by
\begin{eqnarray}\label{6.3}
\breve S(v,p):=\big\{x\in\R^n\big|\;v\in f(x,p)+N_{C(p)}(x)\big\}.
\end{eqnarray}

The main goal of this section is to derive necessary and sufficient conditions for full stability of solutions to \eqref{VC}, \eqref{6.3} expressed entirely via the {\em initial data} of these variational systems.

Given $(\ox,\op)$ satisfying \eqref{6.2}, recall that the partial {\em Mangasarian-Fromovitz constraint qualification} (MFCQ) with respect to $x$ holds at $(\ox,\op)$ if there is $d\in X$ such that
\begin{eqnarray}\label{mfcq}
\la\nabla_x\ph_i(\ox,\op),d\ra<0\;\mbox{ for }\;i\in I(\ox,\op):=\big\{i\in\{1,\ldots,m\}\big|\;\ph_i(\ox,\op)=0\big\}.
\end{eqnarray}
The {\em Lagrangian function} for the variational system \eqref{VC}, \eqref{6.2} is defined by
$$
L(x,p,\lm):=f(x,p)+\disp\sum_{i=1}^{m}\lm_i\nabla_x\ph_i(x,p)\;\mbox{ with }\;x\in\R^n,\;p\in\R^d,\;\mbox{ and }\;\lm\in\R^m.
$$
Under the validity of the partial MFCQ condition \eqref{mfcq}, it is well known that for all feasible $(x,p)$ around $(\ox,\op)$ we have
\begin{eqnarray}\label{6.4}
\begin{array}{ll}
f(x,p)+N_{C(p)}(x)=\big\{L(x,p,\lm)\big|\;\lm\in N\big(\ph(x,p);\Th\big)\big\}:=\Psi(x,p)\;\mbox{ with }\;\Th:=\R^m_-
\end{array}
\end{eqnarray}
and $\ph=(\ph_1,\ldots,\ph_m):\R^n\times\R^d\to\R^m$. Hence any vector $\ov\in f(\ox,\op)+N_{C(\op)}(\ox)$ satisfies
\begin{equation}\label{6.5}
\ov\in L(\ox,\op,\lm)\;\mbox{ with some }\;\lm\in N_\Theta\big(\ph(\ox,\op)\big)
\end{equation}
and the set of Lagrange multipliers is represented by
\begin{eqnarray}\label{6.6}
\Lm(\ox,\op,\ov):=\big\{\lm\in\R_+^m\big|\;\ov\in L(\ox,\op,\lm),\;\la\lm,\ph(\ox,\op)\ra=0\big\}.
\end{eqnarray}
Moreover, it follows from \cite[Proposition~2.2]{LPR} that MFCQ \eqref{mfcq} implies that BCQ \eqref{bcq} holds for $g(x,p)=\delta_{C(p)}(x)$ and that $g$ is parametrically continuously prox-regular at $(\ox,\op)$ for $\hat v=\ov-f(\ox,\op)$.

Let us further recall the following {\em general strong second-order sufficient condition} (GSSOSC) for the variational condition \eqref{VC} as formulated in \cite{K2}: given $(\ox,\op)$ satisfying \eqref{6.2} and given $\ov\in f(\ox,\op)+N_{C(\op)}(\ox)$, the GSSOSC holds at $(\ox,\op,\ov)$ if for all $\lm\in\Lm(\ox,\op,\ov)$ we have
\begin{eqnarray}\label{6.8}
\la u,\nabla_{x}L(\ox,\op,\lm)u\ra>0\;\mbox{ whenever }\;\la\nabla_x\ph_i(\ox,\op),u\ra=0\;\mbox{ as }\;i\in I_+(\ox,\op,\lm),\;u\ne 0
\end{eqnarray}
with the strict complementarity index set $I_+(\ox,\op,\lm):=\big\{i\in I(\ox,\op)\big|\;\lm_i>0\big\}$. This condition is a slight modification and adaptation to \eqref{VC}, \eqref{6.2} of the strong second-order sufficient condition introduced by Robinson \cite{Ro} for nonlinear programs with ${\mathcal C}^2$-smooth data; cf.\ also Kojima \cite{ko}.

Next we modify for the case of PVC in \eqref{VC}, \eqref{6.2} the uniform second-order sufficient condition introduced recently in \cite{MN} for parametric nonlinear programs.

\begin{Definition} {\bf (general uniform second-order sufficient condition).}\label{gusosc} We say that the {\sc general uniform second-order sufficient condition} {\rm (GUSOSC)} holds at $(\ox,\op)$ satisfying \eqref{6.2} with $\ov\in\Psi(\ox,\op)$ if there are positive numbers $\eta,\ell$ such that
\begin{eqnarray}\label{gus}
\begin{array}{ll}
\disp\la\nabla_{x}L(x,p,\lm)u,u\ra\ge\ell\|u\|^2\;\mbox{ for all }\;(x,p,v)\in\gph\Psi\cap\B_\eta(\ox,\op,\ov),\;\lm\in\Lm(x,p,v),&\\
\la\nabla_x\ph_i(x,p),u\ra=0\;\mbox{ as }\;i\in I_+(x,p,\lm)\;\mbox{ and }\;\la\nabla_x\ph_i(x,p),u\ra\ge 0\;\mbox{ as }\;i\in I(x,p)\setminus I_+(x,p,\lm),
\end{array}
\end{eqnarray}
where the mapping $\Psi$ and the set $\Lm(x,p,v)$ are defined in \eqref{6.4} and \eqref{6.6}, respectively.
\end{Definition}

Similarly to the proof of \cite[Proposition~4.2]{MN} in the nonparametric setting we can check that, under the validity of the partial MFCQ \eqref{mfcq}, the GSSOSC from \eqref{6.8} implies the GUSOSC from Definition~{\rm\ref{gusosc}} at $(\ox,\op)$ with $\ov\in\Psi(\ox,\op)$ by passing to the limit in \eqref{gus}.\vspace*{0.03in}

The last qualification condition needed in this section is the partial {\em constant rank constraint qualification} (CRCQ) formulated as follows; cf.\ \cite{FP,L2}. We say that the {\em partial CRCQ} with respect to $x$ holds at $(\ox,\op)$ feasible to \eqref{6.2} if there is a neighborhood $W$ of $(\ox,\op)$ such that for any subset $J$ of $I(\ox,\op)$ the gradient family $\big\{\nabla_x\ph_i(x,p)\big|\;i\in J\big\}$ has the same rank in $W$. It occurs that the simultaneous fulfillment of the partial MFCQ and CRCQ ensures the Lipschitz-like property of the graphical mapping $K$ from \eqref{K} crucial for the (Lipschitzian) full stability results in Section~4.

\begin{Proposition}{\bf (graphical Lipschitz-like property under partial MFCQ and CRCQ).}\label{cr} Assume that both partial MFCQ and CRCQ conditions hold at the point $(\ox,\op)$ feasible to \eqref{6.2}. Then, given any vector $\ov$ satisfying \eqref{6.5}, the Lipschitz-like property in \eqref{K} holds around $(\op,\ox,\hat v)$ with $g(x,p)=\delta_{C(p)}(x)$ and $\hat v=\ov-f(\ox,\op)$.
\end{Proposition}
{\bf Proof.} Follows directly from \cite[Proposition~5.2]{MN} for the case of constant cost functions.\endproof

Now we are ready to derive the main result of this section.

\begin{Theorem}{\bf (characterization of full stability for PVC under partial MFCQ and CRCQ).}\label{thm6.3} Take $(\ox,\ov,\op)\in\R^n\times\R^n\times\R^d$ with $\ov\in\breve S(\ox,\op)$ from \eqref{6.3} and suppose that both partial MFCQ and CRCQ conditions hold at $(\ox,\op)$. Then the following assertions are equivalent:

{\bf (i)} $\ox$ is fully stable solution of the PVC in \eqref{VC} and \eqref{6.2} corresponding to $(\ov,\op)$.

{\bf(ii)} The GUSOSC from Definition~{\rm\ref{gusosc}} holds at $(\ox,\op)$ with $\ov\in\Psi(\ox,\op)$.\\
Consequently, the validity of GUSOSC at $(\ox,\op)$ with $\ov\in\Psi(\ox,\op)$ ensures that the solution map $\breve S$ admits a single-valued and Lipschitz continuous localization $\vt$ around $(\ov,\op,\ox)$.
\end{Theorem}
{\bf Proof.} It suffices to show by Theorem~\ref{coro2} that the imposed GUSOSC is equivalent to the second-order subdifferential condition \eqref{4.8}. This can be done via calculating the term $\la z,w\ra$ in \eqref{4.8} by using the formula for $\Hat D^*\partial g_p$ established in \cite[Theorem~6]{HKO} under the validity of MFCQ and CRCQ. The proof is similar to the one given in \cite[Theorem~5.3]{MN}, and so we omit details. \endproof

Another sufficient condition for the existence of a single-valued and Lipschitz continuous localization of the solution map $\breve S$ around $(\ov,\op,\ox)$ was obtained by Facchinei and Pang \cite{FP} under the name of the {\em strong coherent orientation condition} (SCOC) from \cite[Definition~5.4.11]{FP}, with imposing both MFCQ and CRCQ while assuming in addition that all the functions $\ph_i$ in \eqref{6.2} are convex. The latter convex assumption has been dismissed in the more recent paper by Lu \cite{L2}. Observe that both developments in \cite{FP,L2} rely on topological degree theory the application of which to sensitivity analysis in optimization was initiated by Kojima \cite{ko}.

It is worth mentioning the previous result in this direction by Kyparisis \cite{K2} who proved, based mainly on the implicit function technique by Robinson \cite{R}, the local single-valuedness and continuity (while not Lipschitz continuity) of the solution map $\breve S$ for convex PVC \eqref{VC} (i.e., for quasi-variational inequalities) under the so-called ``general modified strong second-order condition" that is stronger than GSSOSC \eqref{6.8}, which in turn is known to be stronger than SCOC.

The next example demonstrates that our new GUSOSC from Definition~\ref{gusosc} is {\em strictly weaker} than GSSOSC and is {\em not implied} by SCOC even in the case of constraint functions linear with respect to the decision $x\in\R^3$ and parameter $p\in\R^2$ variables as well as of cost functions linear in $p$ and quadratic in $x$ under the validity of both MFCQ and CRCQ conditions. This shows that Theorem~\ref{thm6.3}, which completely characterizes the {\em stronger} property of full stability of solutions to nonconvex PVC derived via advanced techniques of variational analysis and generalized differentiations, provides new sufficient conditions for the local single-valuedness and Lipschitz continuity of the solution map \eqref{6.3} independent of \cite{FP,L2} and significantly extended those in \cite{K2}.

\begin{Example}{\bf (improving sufficient conditions for single-valuedness and Lipschitz continuity of PVC under partial MFCQ and CRCQ).}\label{ex}
{\rm Consider the PVC in \eqref{VC}, \eqref{6.2} with
\begin{eqnarray}\label{5.30}
\left\{\begin{array}{ll}
f(x,p)=\nabla\ph_0(x,p)\;\mbox{ for }\;\ph_0(x,p):=x_3+\Big(\frac{1}{4}+p_2\Big)x_1+p_1x_2+x^2_3-x_1x_2,\\
\ph_1(x,p):=x_1-x_3-p_1\le 0,\\
\ph_2(x,p):=-x_1-x_3+p_1\le 0,\\
\ph_3(x,p):=x_2-x_3-p_2\le 0,\\
\ph_4(x,p):=-x_2-x_3+p_2\le 0,\\
x=(x_1,x_2,x_3)\in\R^3,\;p=(p_1,p_2)\in\R^2.
\end{array}\right.
\end{eqnarray}
We get that both partial MFCQ and CRCQ hold at $(\ox,\op)$ with $\ox=(0,0,0)$, $\op=(0,0)$. It follows from the arguments in \cite[Example 6.4]{MN} that our GUSOSC holds in this setting. Thus we have by Theorem~\ref{thm6.3} that $\ox$ is fully stable in \eqref{5.30} and consequently the solution map $\breve S$ admits a single-valued Lipschitz continuous localization $\vt$ around $(\ov,\op,\ox)$ with $\ov=(0,0,0)$. We also deduce from \cite[Example~6.4]{MN} that GSSOSC \eqref{6.8} fails when $\bar\lm:=(\frac{3}{8}, \frac{5}{8},0,0)$.

Let us now check that the aforementioned SCOC does not hold here. Indeed, note that this vector $\bar\lm$ is an extreme point of the set $\Lambda(\ox,\op,\ov)$, which can be calculated directly by \eqref{6.6} as
$$
\Lambda(\ox,\op,\ov)=\Big\{\Big(\frac{3}{8}-\al,\frac{5}{8}-\al,\al,\al\Big)\Big|\;0\le\al\le\frac{3}{8}\Big\}.
$$
Observe further that the gradient vectors $\nabla_x\ph_1(\ox,\op),\nabla_x\ph_2(\ox,\op)$ are linearly independent in $\R^3$ and that the determinant of the matrix
\[\begin{pmatrix}
\nabla_x L(\ox,\op,\bar\lm)&\nabla\ph_1(\ox)^T&\nabla\ph_2(\ox)^T\\
-\nabla\ph_1(\ox)&0&0\\
-\nabla\ph_2(\ox)&0&0
\end{pmatrix}=\begin{pmatrix}
0&-1&0&1&-1\\
-1&0&0&0&0\\
0&0&2&-1&-1\\
-1&0&1&0&0\\
1&0&1&0&0
\end{pmatrix}
\]
is equal to zero. It shows the violation of SCOC from \cite[Definition~5.4.11]{FP} in this example.}
\end{Example}

Finally in this section, we discuss the role of the {\em pointwise} second-order condition GSSOSC in the study of full stability of PVC from \eqref{VC}, \eqref{6.2}. The following consequence of Theorem~\ref{thm6.3} describes the situation under the two first-order constraint qualifications considered above as well as under the partial {\em linear independence constraint qualification} (LICQ) at $(\ox,\op)$:
\[
\mbox{the gradients}\;\nabla_x\ph_i(\ox,\op)\;\mbox{ for }\;i\in I(\ox,\op)\;\mbox{ are linearly independent},
\]
which clearly implies the partial MFCQ and CRCQ at the corresponding point. Note that both assertions of the corollary below are new and rather surprising by taking into account the classical nature of the first-order and second-order conditions used and the novelty of the full stability notion for general PVC under consideration. On the other hand, we have recently obtained in \cite{MN} the prototypes of these results for parametric nonlinear programs.

\begin{Corollary} {\bf (full stability of PVC via GSSOSC).}\label{gssosc} Let $(\ox,\ov,\op)\in\R^n\times\R^n\times\R^d$ be such that $\ox\in\breve S(\ov,\op)$ in \eqref{6.3}. The following assertions hold:

{\bf (i)} If both partial MFCQ and CRCQ are satisfied at $(\ox,\op)$, then the validity of GSSOSC \eqref{6.8} ensures that $\ox$ is a fully stable solution of PVC corresponding to $(\ov,\op)$.

{\bf (ii)} If the partial LICQ holds at $(\ox,\op)$, then GSSOSC is necessary and sufficient for the full stability of the corresponding solution $\ox$ in {\bf(i)}.
\end{Corollary}
{\bf Proof.} Since GSSOSC is stronger than GUSOSC as discussed after Definition~\ref{gusosc}, assertion {\bf (i)} follows directly from Theorem~\ref{thm6.3}. To verify the necessary of GSSOSC for full stability in assertion {\bf (ii)} under the partial LICQ, it suffices to show that the pointwise second-order subdifferential condition \eqref{4.43} from Theorem~\ref{Lips} is equivalent to GSSOSC. The proof of this fact follows the lines in the proof of \cite[Theorem~6.6]{mrs}, where the inner product $\la z,w\ra$ with $(z,w)\in(D^*\partial_x g)(\ox,\op,\hat v)(w)$ and $\hat v=\ov-f(\ox,\op)$ is explicitly calculated.\endproof\\
{\bf Acknowledgements.} The authors express their deep gratitude to Heinz Bauschke and Shawn Wang for insightful discussions on local maximal monotone operators.

\end{document}